\documentclass[12pt]{amsart}
\usepackage{amsthm}
\input{epsf.tex}
\newtheorem{theorem}{Theorem}[section]
\newtheorem{proposition}[theorem]{Proposition}
\newtheorem{lemma}[theorem]{Lemma}

\newtheorem{define}[theorem]{Definition}

\def\Empty{}

\catcode`\@=11

\def\section{\@startsection {section}{1}{\z@}{-3.5ex plus -1ex minus 
-.2ex}{2.3ex plus .2ex}{\large\bf}}

\def\fnum@figure{{\small Figure \thefigure}}
\def\fakefigure{\def\@captype{figure}}

\long\def\@makecaption#1#2{
    \vskip 10pt 
    \def\FCap{#2} \def\NoCap{\ignorespaces}
    \ifx \FCap\NoCap
       \setbox\@tempboxa\hbox{#1}  
      \else
       \setbox\@tempboxa\hbox{#1: \small \it #2}
    \fi
    \ifdim \wd\@tempboxa >\hsize   
        \unhbox\@tempboxa\par      
      \else                        
        \hbox to\hsize{\hfil\box\@tempboxa\hfil}  
    \fi}

\catcode`\@=12

\def\oplabel#1{
  \def\OpArg{#1} \ifx \OpArg\Empty {} \else
  	\label{#1}
  \fi}

\newlength{\saveu}

\input{serg.tex}
\def\leftceiling{{\lceil}}
\def\rightceiling{{\rceil}}
\def\Reals{{\bf R}}
\def\reals{\Reals}
\def\Integers{{\bf Z}}
\def\integers{\Integers}

\def\Hyperbolic{{\bf H}}
\def\hyperbolic{\Hyperbolic}
\def\Euclidean{{\bf E}}

\def\Complex{{\bf C}}
\def\complex{\Complex}
\def\union{\cup}
\def\Union{\bigcup}
\def\intersect{\cap}

\def\Intersection{\bigcap}
\def\inverse{{-1}}
\def\boundary{\partial}

\def\bdy{\boundary}

\def\composed{\circ}

\def\restrictedto{\bigm|}

\def\cross{\times}

\def\infinity{\infty}

\def\Sum{\sum}

\def\from{\colon}
\def\homeo{\approx}

\def\func#1{\operatorname{#1}}

\def\cl{\func{cl}}
\def\closure{\cl}

\def\Isom{\func{Isom}}
\def\interior{\func{int}}

\def\Length{\func{Length}}

\def\centeredepsfbox#1{\centerline{\epsfbox{#1}}}

\begin{document}

\title{Quasigeodesic Flows in Hyperbolic Three-manifolds}
\author{S\'{e}rgio Fenley}
\address{Fenley:
Department of Mathematics \\
Univ. of California \\
Berkeley, CA 94720}
\author{Lee Mosher}
\address{Mosher:
Dept. of Math and Computer Science\\
Rutgers University\\
Newark, NJ }
\thanks{Research at MSRI is supported in part by NSF grant DMS-9022140.
Both authors were supported by additional NSF grants.}

\begin{abstract}
Any closed, oriented,  hyperbolic
$3$-manifold with nontrivial second homology has many quasigeodesic
flows, where quasigeodesic means that flow lines are uniformly efficient
in measuring  distance in relative homotopy classes. The flows are
pseudo-Anosov flows which are  almost transverse to finite depth
foliations in the manifold. The main tool is the use of a sutured
manifold hierarchy which has good geometric properties.
\end{abstract}

\maketitle

\section*{Introduction}

In this article we prove that any closed, oriented,  hyperbolic
$3$-manifold with nontrivial second homology has many quasigeodesic
flows, where quasigeodesic means that flow lines are uniformly efficient
in measuring  distance in relative homotopy classes. The flows are
pseudo-Anosov flows which are  almost transverse to finite depth
foliations in the manifold. The main tool is the use of a sutured
manifold hierarchy which has good geometric properties.

The best metric property a flow can have is that all its flow lines are
minimal geodesics in their relative homotopy classes, which amounts to
being minimal geodesics when lifted to the  universal cover. Suspensions
of Anosov diffeomorphisms of the torus and geodesic flows on the unit
tangent bundle of surfaces of constant negative curvature have this
property. Even for these examples one has to choose an appropriate
metric to get the minimal property. If the metric is changed the flow
lines are only quasigeodesics: when lifted to the universal cover,
length along flow lines is a bounded multiplicative distortion of length
in the manifold plus an additive constant. The concept of quasigeodesic
has the advantage of being independent of the metric in the manifold. We
say that a flow is quasigeodesic if all flow lines are quasigeodesics.

Our main interest is in hyperbolic manifolds. In these manifolds, a
quasigeodesic in the universal cover is a bounded distance from some
minimal geodesic \cite{Th1}. This, among other reasons, makes
quasigeodesics extremely useful in studying hyperbolic manifolds
\cite{Th1,Th2, Mor,Ca}.

A natural question to ask is: how common are quasigeodesic flows? Notice
that Zeghib \cite{Ze} proved that there cannot exist a continuous
foliation by geodesics in a hyperbolic $3$-manifold. On the other hand,
in their seminal work \cite{Ca-Th}, Cannon and Thurston showed that if a
hyperbolic 3-manifold $M$ fibers over the circle, then $M$ has a
quasigeodesic flow which is transverse to the fibers. Afterwards Zeghib
\cite{Ze} gave a quick and elementary proof that for any compact manifold
$M$ which fibers over a circle, any flow transverse to the fibration is
quasigeodesic. Mosher \cite{Mo3} produced examples of quasigeodesic 
flows transverse to a class of depth one foliations in hyperbolic
$3$-manifolds.

The flow constructed in \cite{Ca-Th} is the suspension of a
pseudo-Anosov homeomorphism of the fiber. Hence it is a {\em
pseudo-Anosov} flow, that is, it has stable and unstable foliations in
the same way as Anosov flows do, except that one allows $p$-prong
singularities along finitely many closed orbits. The quasigeodesic and
pseudo-Anosov properties are used in an essential way in
Cannon-Thurston's proof that lifts of fibers extend continuously to the
sphere at infinity, providing examples of sphere filling curves. The
flows constructed in \cite{Mo3} are also pseudo-Anosov. Mosher showed
that quasigeodesic pseudo-Anosov flows on hyperbolic manifolds can be
used to compute the Thurston norm \cite{Mo1,Mo2}.

The quasigeodesic property for Anosov flows in hyperbolic $3$-manifolds
has also been extensively studied by Fenley who showed that there are
many examples which are not quasigeodesic \cite{Fe1}. In addition the
quasigeodesic property for Anosov flows is related to the topology of
the stable and unstable foliations in the universal cover \cite{Fe2}
and implies that limit sets of leaves of these foliations are Sierpinski
curves \cite{Fe3}. 

The main goal of this paper is to show that quasigeodesic flows are
quite common. If $M^3$ is closed, oriented, 
irreducible with $H_2(M) \not = 0$ and
if $\zeta \not = 0$ in $H_2(M)$, then Gabai \cite{Ga1} constructed a
taut, finite depth foliation $\fol$ whose set of compact leaves
represents $\zeta$. Given such $\fol$ in a hyperbolic $3$-manifold,
Mosher \cite{Mo4} constructed pseudo-Anosov flows which are almost
transverse to $\fol$. Almost transverse means that it will be transverse
to $\fol$ after an appropriate blow up of a finite collection of closed
orbits (see detailed definition in section $4$).

\medskip
\noindent {\bf Main theorem} \ {\em Let $M$ be a closed, oriented,
hyperbolic
$3$-manifold with non zero  second betti number
and let $\zeta$ a nonzero homology
class in $H_2(M)$. 
Let $\fol$ be a taut, finite depth foliation whose compact
leaves represent $\zeta$, and $\Phi$ a pseudo-Anosov flow which is
almost transverse to $\fol$. Then $\Phi$ is a quasigeodesic flow.}
\medskip

If the compact leaves of $\fol$ are fibers, the main theorem follows from
the above mentioned result of Zeghib by showing that every flow line of
$\Phi$ hits the compact leaves. If some compact leaf is not a fiber the
theorem easily follows from a more general result:

\medskip
\noindent {\bf Theorem A} \
{\em Let $M$ be a closed, oriented,  hyperbolic $3$-manifold with
a taut, finite depth foliation $\fol$, so that some 
compact leaf is not a fiber of $M$ over the circle. Let $\Phi$
be a flow transverse to $\fol$ and let $\wwp$ be the
lifted flow to the universal cover.
Then $\Phi$ is quasigeodesic if and only if
$\wwp$ has Hausdorff orbit space.
}
\medskip

The orbit space of a flow is the quotient space obtained by collapsing
each flow line to a point. 
We remark that the only if part of theorem A is straightforward.  One
might also ask whether the condition in theorem $A$ is non void, that
is, if there are flows transverse to Reebless finite depth foliations
for which the orbit space of $\wwp$ is not Hausdorff. Indeed this is
possible. Examples are easy to construct containing a flow invariant
annulus $Z$ in $M$, so that flow lines induce in $Z$ a $2$-dimensional
Reeb foliation. In that situation clearly the orbit space of $\wwp$ is
not Hausdorff.

When $\fol$ has a nonfiber compact leaf the main theorem follows from
theorem A via two remarks. (1) The important fact is that for the
original pseudo-Anosov flow $\Phi$, the covering flow
$\wwp$ has Hausdorff orbit space. This is where the pseudo-Anosov
dynamics plays an essential role. (2) This implies that the blown up
flow also has Hausdorff orbit space, and since it is transverse to
$\fol$ it is quasigeodesic by theorem A. As a consequence the original
flow is also quasigeodesic. 

Here are the key ideas in the proof of theorem $A$. 
Some (and hence any) compact leaf of
$\fol$ 
represents a quasi-Fuchsian subgroup. Therefore this leaf lifts to a
quasi-isometrically embedded surface in $\mi$, which has excellent
geometric properties. 
That means that if a 
flow line in $\mi$ keeps intersecting lifts
of compact leaves, these lifts trap the flow line
which then converges to a single
point in $\si$.

Next we proceed to extend this argument to all orbits.
For that, we use sutured manifold hierarchies
and branched surfaces associated to the foliation $\fol$. By general
principles the sutured manifolds in the hierarchy have good geometric
properties, that is, they are quasi-isometrically embedded when lifted
to the universal cover and the cutting surfaces in the hierarchy are also
quasi-isometrically embedded. 
The cutting surfaces play the role of compact surfaces in
the appropriate sutured manifold in the hierarchy.
Using an induction argument with the sutured manifold hierarchy,
we can show that for any point in $\mi$, its
flow line in forward time converges to a unique point in the sphere at
infinity $\si$, and likewise for the negative direction.  
This is first shown in the compactified universal cover of 
the appropriate sutured manifold in the hierarchy and then
derived in $\hhh \cup \si$ by way of the good geometric properties
of the sutured manifolds.
The existence of a unique limit point of flow lines
is a much
weaker property than being quasigeodesic: for example horocycles have
this property but they are not quasigeodesic.

We also show that the 
limit point map of flow lines
is continuous and that forward and backward limit points
in each orbit are distinct. Quasigeodesics satisfy all of these
properties. These results only keep track of the asymptotic behavior of
a flow line, but a priori do not determine the rough location of the
flow line --- which must be the case for quasigeodesics.  For
flows, however, these properties are indeed sufficient to ensure
quasigeodesic behavior as proved by:

\medskip
\noindent {\bf Theorem B} \
{\em 
Let $\Phi$ be a flow in $M^3$ closed hyperbolic.
Suppose that:

(a) Each half orbit of $\wwp$ has a unique limit point in $\si$,

(b) For a given orbit, the forward and backward limit points
are distinct,

(c) The forward and backward limit point maps are continuous.

\noindent
Then $\Phi$ is quasigeodesic.
}
\medskip

The paper is organized as follows.
In section $1$ we prove a generalization of theorem $B$, 
which applies to closed
invariant sets of $\Phi$. This is used in the inductive step of
the proof of theorem $A$.
In section $2$ we study sutured manifold hierarchies
adapted to finite depth 
foliations and prove the needed geometric properties
of the sutured manifolds.  
In section $3$ we prove that conditions (a,b,c) of theorem
$B$ hold for flows satisfying the hypothesis of theorem $A$,
and as a consequence derive theorem $B$.
In the final section we study pseudo-Anosov flows and prove
the main theorem.

We thank Dick Canary who pointed out to us the needed quasi-isometric
properties of finitely generated Kleinian groups and also for many
useful suggestions.

\vskip .1in
\section{From continuous extension to quasigeodesic behavior}
\vskip .1in

Here are the basic definitions we need concerning quasi-isometries
and quasigeodesics. A good source for foundational material
on quasi-isometries is \cite{BKS}, especially chapter 10 by
Ghys and De la Harpe on ``Infinite groups as geometric objects (after
Gromov)'' and chapter 11 by Cannon on ``The theory of negatively
curved groups and spaces''. 

\begin{define}{(quasi-isometries and quasigeodesics)}{}
Given metric spaces $(X,d)$ and $(Y,d')$, a map $f \from X \to Y$ is a
quasi-isometric embedding if there is $k \geq 1$ so that for any $z, w \in
X$, we have

$$ \frac{1}{k} \, d'(f(z),f(w)) - 1 \leq d(z,w) \leq k \, d'(f(z),f(w)) +
k.$$

\noindent
Once a metric is fixed we say that $f$ is a $k$-quasi-isometry. 

The spaces $X,Y$ are quasi-isometric if there is a quasi-isometry $f
\from X \to Y$ and a constant $k'$ such that each point of $Y$ is
within distance $k'$ of the image of $f$.

Given a metric space $X$, a quasigeodesic is a map
$f \from \rrrr \to X$ which is a quasi-isometry. If $f$ is a
$k$-quasi-isometry we say that $f$ is a $k$-quasigeodesic.
\end{define}

A {\em flow} on a manifold $X$ is a continuous action of $\reals$ on $X$,
i.e.\ a continuous map $\psi \from X \cross \reals \to X$ written $(x,t) \to
\psi_t(x)$, such that 
$$\psi_{s+t}(x) = \psi_s(\psi_t(x)) \quad\text{for all}\quad x \in X, \,
s,t \in \reals.
$$
A {\em semiflow} on $X$ is defined similarly, except that the domain of
$\psi$ is a closed subset $D \subset X \cross \reals$ such that for each $x
\in X$, $(x\cross \reals) \intersect D = x \cross J$ for some closed,
connected set $J \subset \reals$ containing $0$, and the above equation holds
whenever the two sides are defined.  All flows and semiflows in this article
are generated by nonzero, continuous vector fields, i.e.\ flow lines are
smooth immersions and the tangent vector field is continuous on $X$. When
$\psi$ is understood, we often write $x \cdot t = \psi_t(x)$, and if $J
\subset \reals$ is an interval we write $x \cdot J = \psi_J(x)$. If $y = x
\cdot t$ let $[x,y] = x \cdot [0,t]$ and $(x,y) = [x,y] - \{ x, y \}$; also
let $\tau(x,y) = |t|$.

\begin{define}{(uniformly quasigeodesic flows)}{}
Let $\Phi$ be a flow on a compact manifold $M$. Let $\widetilde \Phi$ be the
lifted flow on the universal cover $\widetilde M$. We say that $\Phi$ is
uniformly quasigeodesic if there exists a constant $k \ge 1$ such that for
each $x\in \widetilde M$, the map $t \to x \cdot t$ is a
$k$-quasigeodesic. 
\end{define}

Since $M$ is compact this property is independent of the parameterization of
$\Phi$. Henceforth flowlines are always parameterized by arc length. Given
this, a flow line in $\widetilde M$ is $k$-quasigeodesic if for any $x,y$ in
the flow line, $\tau(x,y) \le k \, d(x, y) + k$. We sometimes say that flow
lines of $\Phi$ are quasigeodesic.

We first show that there is a weak form of quasigeodesic
behavior which follows from purely topological properties of
the orbit space of $\wwp$.
Given $\epsilon, T > 0$ an $\epsilon, T$ cycle of $\Phi$ is a closed
loop in $M$ obtained from an orbit segment of length $> T$ with endpoints
less than $\epsilon$ apart, closed up by an arc of length $< \epsilon$.
The following lemma does not assume that the manifold is
hyperbolic.
This lemma is not logically necessary for proving the main theorem
and theorems $A, B$;
but it describes additional hypotheses (which are quite common)
under which the proofs of these theorems can be simplified.

\begin{lemma}{(weak quasigeodesics)}{}
Let $Y$ be a closed invariant set of a flow $\Phi$ in a closed manifold $M$.
Let $\widetilde Y$ be the lift of $Y$ to $\mi$. Suppose that 

(a) $\wwp |_{\widetilde Y}$ has Hausdorff orbit space,

(b) There are $\epsilon, T > 0$ so that any
$\epsilon, T$ cycle of $\Phi |_Y$ is not null homotopic in $M$.

\noindent
Then for any $b > 0$ there is $c_b > 0$ (depending only on $b$)
so that if $x, y$ are in an orbit $\gamma$ of $\wwp$ and
$\tau(x,y) > c$ then $d(x,y) > b$.
\label{weakqg}
\end{lemma}

\begin{proof}{}
Otherwise there is $b > 0$ and $x_i, y_i \in \widetilde Y$,
with $x_i, y_i$ in the same orbit $\gamma_i$ of $\wwp$ and
$\tau(x_i,y_i) > i$ but $d(x_i,y_i) < b$.
Since $M$ is compact, then up to covering translations
and taking subsequences we may assume that $x_i \rightarrow x$
and $y_i \rightarrow y$.
Notice that $x, y \in \widetilde Y$ since $\widetilde Y$ is closed.

By hypothesis (a) it follows that $x$ and $y$ are in the same
orbit of $\wwp$ so $y = x \cdot t$ for some $t \in \rrrr$.
By the local product structure of the flow along compact
orbit segments, there are $t_i \rightarrow t$ so
that $z_i = x_i \cdot t_i \rightarrow y$.
If $i$ is big enough then $\tau(z_i,y_i) > T$
and $d(z_i,y_i) < \epsilon$, thereby producing 
$\epsilon, T$-cycles
of $\wwp$. These project to null homotopic
$\epsilon, T$-cycles of $\Phi$, contradiction.
\end{proof}

\noindent
{\bf Remarks:} (1) Notice in particular that
condition (b) implies that orbits of $\wwp$ are
never periodic and are properly embedded in $\mi$.
When (b) holds we say that $\Phi |_Y$ satisfies the
$\epsilon, T$-cycles condition.

(2) Quasigeodesic behavior is the additional property that
$c_b$ is bounded by an affine function of $b$.
In general conditions (a),(b) are not sufficient to ensure
quasigeodesic behavior. For instance Anosov flows
always satisfy these conditions \cite{Fe1}, but there are
many examples of Anosov flows in hyperbolic manifolds which
are not quasigeodesic \cite{Fe1}.

(3) When $Y=M$ 
conditions (a),(b) together are equivalent to the orbit space $\oo$ of
$\wwp$ being homeomorphic to either the plane $\rrrr^2$
or the sphere ${\bf S}^2$.
(a) implies that $\oo$ is a 2-dimensional manifold,  and it 
is Hausdorff. Since $\oo$ is simply connected and
has no boundary it is either ${\bf S}^2$ or $\rrrr^2$.
This means that the flow $\wwp$ 
is topologically a {\em product} flow in $\mi$.
Lemma \ref{weakqg} means that topological product flows
always satisfy a weak quasigeodesic property.
This is reminiscent of the same situation for codimension
one foliations which was studied in \cite{Fe4}.
Finally we remark that our main interest is in
closed hyperbolic
manifolds, where $\pi_2(M)$ is trivial.
In that case conditions (a),(b) of lemma \ref{weakqg} are
equivalent to $\wwp$ having orbit space homeomorphic to $\rrrr^2$.

\smallskip
For flows in hyperbolic $3$-manifolds, 
we now develop a method to  upgrade
information about asymptotic behavior of flow lines
of $\wwp$ into metric efficiency of flow lines.
This will be the key tool to prove uniform quasigeodesic
behavior for a large class of flows in hyperbolic $3$-manifolds.

\begin{theorem}{}{}
Let $Y$ be a closed invariant set of a non singular flow $\Phi$ 
in $M^3$ closed hyperbolic. Suppose that



(a) half orbits converge: for any $x \in \widetilde Y$
each of the two rays of $x \cdot \rrrr$ accumulate in a single 
point of $\si$, that is,
the following limits exist:

$$\lim_{t \rightarrow +\infty} x \cdot t = \eta_+(x) \in \si
\ \ {\rm and} \ \ 
\lim_{t \rightarrow -\infty} x \cdot t = \eta_-(x) \in \si.$$

(b) For each $x \in \widetilde Y$, $\eta_+(x) \not = \eta_-(x)$.

(c) The maps $\eta_+, \eta_-: \widetilde Y \rightarrow \si$
are continuous.

\noindent
Then the orbits  in $\wwp |_{\widetilde Y}$ are uniformly
quasigeodesics: there is $k > 0$ so that for any 
orbit $\gamma$ of $\wwp$ in $\widetilde Y$ and for
any $x, y \in \gamma$, 
$\tau(x,y) < k d(x,y) + k$.
\label{upgrade}
\end{theorem}

Notice that conditions (a), (b) and (c) are necessary to
get quasigeodesic behavior. Conditions (a) and (b) follow
directly from the fact that single flow lines are quasigeodesics
\cite{Th1,Gr,Gh-Ha}. Condition (c) is not true for an 
arbitrary collection of quasigeodesics, but holds for 
uniformly quasigeodesic closed invariant sets of flows.
Since we will use this last fact in the proof of theorem $A$, we
provide a proof in section \ref{qg}.

\begin{proof}{}
Since $\eta_+(x) \not = \eta_-(x)$ for  any $x \in \widetilde Y$,
let $g_x$ be the unique geodesic in $\hhh$ with endpoints
$\eta_+(x), \eta_-(x)$  and
let $\bdy g_x \subset \si$ be the ideal points of $g_x$.
Given $a > 0$ let $U_a(g_x) \subset \hhh$ be the neighborhood
of radius $a$ around $g_x$.
We first show:

\centerline{$(*)$ there exists $a > 0$ so that for any $x \in
\widetilde Y$ we have $x \cdot \rrrr \subset U_a(g_x)$.}

\noindent
Otherwise, let $x_i \in \widetilde Y$ with $d(x_i, g_{x_i}) \to +\infinity$.
Up to covering translations and taking a subsequence assume that $x_i \to x$.
Hence $x \in \widetilde Y$ and by (b), (c) we have 
$$\lim \eta_+(x_i) = \eta_+(x) \ne \eta_-(x) = \lim \eta_-(x_i)
$$ 
But $d(x,g_{x_i}) \to +\infinity$,
so up to taking a subsequence we may assume
that there exists $p
\in S^2_\infinity$ such that $\lim (g_{x_i} \cup \partial g_{x_i})  = p$ 
in the Hausdorff topology
on closed subsets of $\hhh \union S^2_\infinity \homeo B^3$. 
Therefore $p = \lim \bdy g_{x_i} = \lim \{ \eta_+(x_i), \eta_-(x_i) \}$,
contradicting that $\lim \eta_+(x_i) \ne \lim \eta_-(x_i)$.

\smallskip
We now assume that $\wwp | _{\widetilde Y}$ is not uniformly 
quasigeodesic and derive a contradiction.
There are two steps in the argument:

\begin{enumerate} 

\item[] Step(1) Using $(*)$ we show that there are  $x_i, y_i$ in the
same orbit of $\widetilde Y$ with  $d(x_i,y_i)$ bounded but
$\tau(x_i,y_i) \rightarrow +\infty$. 

\end{enumerate}

\noindent
If in addition one knows that $\Phi |_Y$ satisfies the
$\epsilon,T$-cycles condition and
$\wwp | _{\widetilde Y}$ has 
Hausdorff orbit space, then the conclusion of step
(1) is disallowed by lemma \ref{weakqg}, finishing the proof.
However, these additional assumptions
are not necessary because:

\begin{enumerate}
\item[] Step (2) The conclusion of step (1) is disallowed by conditions
(a), (b) and (c).
\end{enumerate}

We need the following definitions: given $w \in \hhh$ and an oriented
geodesic $g \subset \hhh$, let $P(w,g)$ be the hyperplane of $\hhh$ 
containing $w$ and perpendicular to $g$. Let $F(w,g)$ be the component of
$\hhh - P(w,x)$ containing the positive endpoint of $g$ in its closure and
let $B(w,x)$ be the other component. For $w \in g$ and $b > 0$ let $w+b$
be the point of $g$ with $d(w, w+b) = b$ and $w+b$ in the ray from $w$ to
the positive endpoint of $g$.

For any $x \in \widetilde Y$ we define $\rho_x \from \hhh \to g_x$ to be the
closest point projection onto $g_x$. Then for any $x, y$ in the same orbit of
$\wwp$ in $\widetilde Y$,
$$d(\rho_x (x), \rho_x (y)) \le d(x,y) \le 
d(\rho _x (x), \rho _x (y)) + 2a \eqno(1)
$$ 
because $x,y \in U_a(g_x)$.

To prove step (1), because we have assumed that orbits are not
uniformly quasigeodesic, for any
$C > 0$ there is an orbit segment $\gamma = [x,y]$ in $\widetilde Y$
which satisfies
$\Length(\gamma) / d(x,y) > 2C$ and $\Length(\gamma) > C$. 
Assume first that $d(x,y) \ge 1 + 2a$.
Hence $d(\rho _x (x), \rho _x (y)) \ge 1$ by (1).
Also $d(x,y) \ge
d(\rho _x (x), \rho _x (y))$, so

\begin{eqnarray*}
{\frac{\Length(\gamma)}{d(\rho _x (x), \rho _x (y))}} \ 
	\geq \ \frac{\Length(\gamma)}{d(x,y)} \  > \ 2C \
\ge \ C + \frac{C}{d(\rho _x (x), \rho _x (y))} \\
\end{eqnarray*}

\noindent
and therefore

\begin{eqnarray*}
\frac{\Length(\gamma)}{C} \ >  \ d(\rho _x (x),\rho _x (y)) + 1 \
  > \ \leftceiling d(\rho _x (x),\rho _x (y)) \rightceiling \\
\end{eqnarray*}

\noindent
where $\leftceiling x \rightceiling$ denotes the ``ceiling'' function, 
that is, the
least integer $\ge x$. Set $n_0 = \leftceiling d(\rho _x (x),\rho _x (y))
\rightceiling$, so
$\Length(\gamma) > n_0C$. Also, 
$(n_0-1) < d(\rho _x (x),\rho _x (y)) \le n_0$,
so we can find consecutive points 
$\rho _x (x) = z_0, z_1, \ldots, z_{n_0} = \rho _x (y)$
on $g_x$ such that $d(z_{n-1},z_n) = 1$ if $1 \le n < n_0$ and
$0 < d(z_{n_0-1},z_{n_0}) \le 1$. 
Let $x=x_0$ and for $1 \le n \le n_0$ let $x_n$
be the last point on $\gamma = [x,y]$ such that $\rho _x (x_n) = z_n$, so
the points $x_0,x_1,\ldots,x_{n_0}$ are consecutive on $\gamma$ and they
partition $\gamma$ into subsegments $\gamma = \gamma_1 * \cdots * \gamma_{n_0}$
with $\gamma_n = [x_{n-1},x_n]$. Since 
$$\Sum_{n=1}^{n_0} \Length(\gamma_n) =
\Length(\gamma) > n_0C
$$ 
then for some $n$ we have $\Length(\gamma_n) >
C$.

This shows that regardless of whether $d(x,y) \geq 1 + 2a$ or not,
we produce a pair
$x', y'$ so that
$d(x',y') < 1 + 2a$ but $\tau(x',y') > C$.
This finishes the proof of step 1.

We now prove step 2.
Let then $x_i, y_i$ with $d(x_i,y_i) < 1 + 2a$ but
$\tau(x_i,y_i) \rightarrow +\infty$.
Without loss of generality assume that $y_i = x_i \cdot t_i$
with $t_i > 0$.

\smallskip
\noindent
{\underline {Case 1}} --- 
The sequence of intervals 
$[x_i,y_i]$ has a subsequence with bounded diameter, which
we may assume is the original sequence.

Let $v_i$ be the middle point in $[x_i,y_i]$ (with respect to
the parametrization) and up
to subsequence and covering translations assume that $v_i
\rightarrow v \in \widetilde Y$.
For any $t \in \rrrr$, $v_i \cdot t \rightarrow v \cdot t$
and for $i$ big enough $v_i \cdot t \in [x_i,y_i]$.
Since $d(v_i \cdot t_i, v_i)$ is bounded, this shows that
$d(v \cdot t,v)$ is also bounded.
It follows that $v \cdot \rrrr$ is contained in a bounded set
in $\hhh$, hence it accumulates in $\hhh$ contradicting
condition (a) of the theorem.

\smallskip
\noindent
{\underline {Case 2}} --- 
There exist $x_i,y_i \in \widetilde Y$ such that $d(x_i,y_i) < 1 + 2a$
and diameter$([x_i,y_i]) \rightarrow +\infty$. 

Let $g_i = g_{x_i}$ and let
$\rho_i \from \hhh \to g_i$ be the closest point projection. Choose $w_i \in
[x_i,y_i]$ so that

$$d(\rho_i (x_i), \rho_i (w_i))
= \max _{w \in (x_i, y_i)} d(\rho_i(x_i), \rho_i(w)).$$ 

\noindent
Since diam$[x_i,y_i] \rightarrow +\infty$ and $d(x,y) < 1 + 2a$, then
$d(w_i,x_i) >> d(x_i,y_i)$. Assume first that $w_i \in F(x_i, g_i)$, hence
$w_i \in F(y_i, g_i)$. If $x_i \in B(y_i,g_i)$, then as
$P(y_i,g_i)$ separates $x_i$ from $w_i$, we can find $z_i,v_i
\in [x_i,y_i]$ so that:
$$
w_i \in (z_i,v_i), \ \
\rho_i(z_i) = \rho_i(v_i) = \rho_i(y_i)
\ \ {\rm and } \ \ \rho_i(w) \not = \rho_i(y_i), \ 
\forall w \in (z_i,v_i).
$$
\noindent
Similarly if $x_i \in F(y_i,g_i) \cup P(y_i,g_i)$, we can find $z_i, v_i$
satisfying the conditions above except that $\rho_i(z_i) = \rho_i(v_i) =
\rho_i(x_i)$. We obtain a similar statement if $w_i \in B(x_i,g_i)$.

%

These arguments show that there are $x_i, y_i$ so that
$\rho_i(x_i) =
\rho_i(y_i)$, 
$\rho_i(w) \not =
\rho_i(x_i)$ for any $w \in (x_i, y_i)$ and 
diam$[x_i,y_i] \rightarrow +\infty$.
If $(x_i,y_i) \subset F(x_i, g_i)$, 
let $w_i$ be defined as above.
Since $\eta_+(x_i)$ is in the closure of $F(\rho_i(x_i) + 1, g_i)$
in $\hhh \cup \si$ (where $\rho_i(x_i) + 1$ is computed in $g_i$)
and $y_i = w_i \cdot r'_i$ with $r'_i > 0$,
there must be a point  $u_i = w_i \cdot s_i$, with
$s_i > 0$ and $\rho_i(u_i) = \rho_i(w_i)$.
In addition we may assume that $(w_i,u_i) \subset B(w_i,g_i)$
and diam$[w_i,u_i] \rightarrow +\infty$.

\begin{figure}
\centeredepsfbox{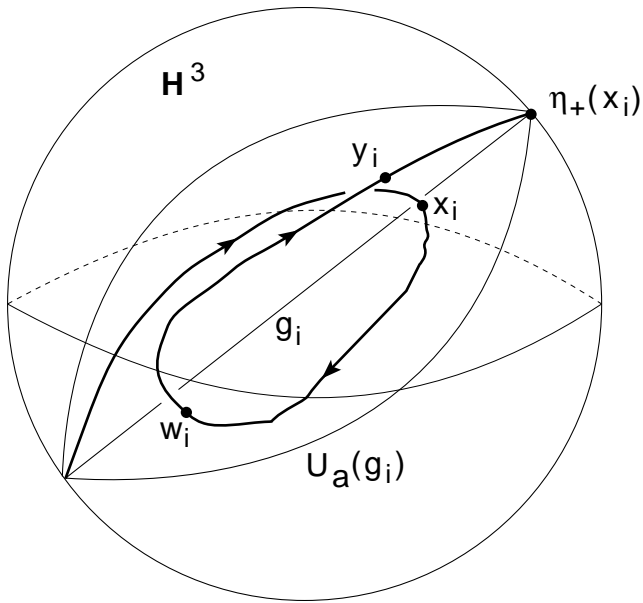}
\caption{Too much backtracking eventually traps orbits in
the negative direction.}
\label{backtrack}
\end{figure}

Since $d(w_i,x_i) \rightarrow +\infty$ these arguments show that in any case
there are $x_i, y_i \in \widetilde Y$, with $\rho_i(x_i) =
\rho_i(y_i)$, $(x_i,y_i) \subset B(x_i,g_i)$ and
diam$[x_i,y_i] \rightarrow +\infty$, see
figure \ref{backtrack}. Up to covering translations and
taking a subsequence we may assume that $x_i \rightarrow x \in \widetilde Y$.
Since $g_i = g_{x_i} \rightarrow g_x$ and $\rho_i(x_i) = \rho_{x_i}(x_i)
\rightarrow \rho_x(x)$, then
$${\overline {P(\rho_{x_i}(x_i) + 2,g_i)}} \ \rightarrow \
{\overline {P(\rho_{x}(x) + 2,g_x)}}$$
\noindent
in the topology of closed sets of $\hhh \cup \si$.
Then
$$\eta_+(x) \in {\overline{F(\rho_{x_i}(x_i) + 1,g_i)}}$$
\noindent 
for $i$ big enough.
Hence there is $t_0 > 0$ so that for $i$ big enough
$x \cdot [t_0,+\infty) \subset F(\rho_{x_i}(x_i),g_i)$. It therefore follows
that for $i$ big enough there is $r_i$ with $r_i < t_0 + 1$ and $x_i \cdot
r_i \in F(\rho_{x_i}(x_i),g_i)$. This contradicts the facts $(x_i,y_i)
\subset B(\rho_{x_i}(x_i), g_i)$ with $y_i = x_i \cdot t_i$ and $t_i
\rightarrow +\infty$. This finally finishes the proof of theorem
\ref{upgrade}.
%
\end{proof}

\noindent
{\bf Remarks:} 

(1) The above proof works verbatim
for $M^n$ closed hyperbolic.



(2) Without much more work, one can also prove this result 
under more general assumptions: 


(i) $M$ is any manifold and $\pi_1(M)$ is word
hyperbolic in the sense of Gromov \cite{Gr}. In that case instead of
$\si$ one uses the appropriate boundary at infinity $\partial
_{\infty} \mi$ --- notice that if $M$ is an oriented, irreducible
3-manifold then $\partial_{\infty} \mi$ is still a $2$-dimensional
sphere \cite{Be-Me}.

(ii) $Y$ is a compact space on which a flow $\Phi$ is defined,
and instead of $Y$ being a subset of $M$ we have instead a continuous
function $f \from Y \to M$. By the pullback construction we have $f
\from \widetilde Y \to \widetilde M$ and a flow $\widetilde \Phi$ on
$\widetilde Y$, and theorem \ref{upgrade} applies.


\newcommand{\R}{{\mathcal R}}
\newcommand{\Fol}{\fol}
\newcommand{\Annuli}{{\mathcal A}}
\newcommand{\bdyv}{\bdy_{v}}
\newcommand{\bdyin}{\bdy_{h}^-}
\newcommand{\bdyout}{\bdy_{h}^+}
\newcommand{\Hull}{{\mathcal H}}
\newcommand{\compact}[1]{\cover #1 \union \bdy_\infinity \cover #1}
\newcommand{\inject}{\hookrightarrow}
\newcommand{\Ker}{\func{Ker}}
\newcommand{\subgroup}{<}
\newcommand{\cover}{\widetilde}
\newcommand{\Rplus}{[0,+\infinity)}
\newcommand{\RPlus}{\Rplus}

\section{Hierarchies}

Finite depth foliations are closely related to sutured manifold
hierarchies \cite{Ga1}. In \cite{Ga3} an ``internal'' version of a
sutured manifold hierarchy was defined, in terms of branched surfaces. We
review this subject here, providing some proofs of known but unpublished
information, and we add some geometric information. For detailed
definitions concerning foliations and laminations  on 3-manifolds see
\cite{Ga1} and \cite{Ga-Oe}.

A $2$-dimensional {\em foliation} of a 3-manifold $M$ is a decomposition
of $M$ into $2$-dimensional manifolds called {\em leaves} which fit
together in a local product structure. A foliation $\fol$ is {\em taut}
if it is  transversely oriented, no leaf is a sphere, and each leaf of
$\fol$  intersects some closed curve in $M$ which is transverse to
$\fol$. The leaves of a taut foliation are $\pi_1$-injective \cite{No}.

A {\em lamination} of $M$ is a foliation of a closed subset of $M$,
covered by charts of the form $D^2 \cross [0,1]$ so that each component of a
leaf intersected with the chart has the form $D^2 \cross t$ for some
$t \in [0,1]$. A lamination $\Lambda$ is {\em essential} if it has no sphere
leaves, no Reeb components, and if $M_\Lambda$ is the metric completion
of $M - \Lambda$, then $M_\Lambda$ is irreducible, boundary
incompressible, and end incompressible; the latter condition means
intuitively that $\bdy M_\Lambda$ has no ``infinite folds''. The leaves
of an essential lamination $\Lambda$ are $\pi_1$-injective in $M$, and
the components of  $M - \Lambda$ are $\pi_1$-injective.

A taut foliation is obviously an essential lamination. It is an exercise
in the results of \cite{Ga-Oe} to show that every sublamination of an
essential lamination is essential.

Given a foliation $\fol$ in a closed manifold $M$  we say that a leaf $L$
of $\fol$ is {\em proper} if $\overline L - L$ is a closed subset  of $M$
and $\fol$ is {\em proper} if  all leaves are proper. In that case Zorn's
lemma implies that $\fol$ has compact leaves, which are then the depth  0
leaves. The {\em depth} is an ordinal $\omega$ defined by induction: a
leaf $L$ is at {\em depth} $\omega$ if $\overline L - L$ is contained in
the union of leaves of depth $< \omega$, but not contained in the union
of leaves $< \omega_1$ for any $\omega_1 < \omega$. A proper foliation
has {\em finite depth} $n$  if $n$ is the maximum of the depth of its
leaves.

A {\em branched surface} in a closed 3-manifold $M$ is a smooth
2-complex $B \subset M$ such that for each $x \in B$, there exists a
smoothly embedded open disc in $B$ containing $x$, and all such discs
are tangent at $x$, therefore determining a well-defined tangent plane
$T_x B$. The set of points where $B$ is not locally a surface is called
the {\em branch locus}. A {\em sector} of $B$ is a complementary
component of the branch locus. In this paper, {\em all branched surfaces
will be transversely oriented}. A flow (or semiflow) is transverse to
$B$ if, for each $x \in B$, the tangent vector to the flow points
to the same side of $T_x B$ as the transverse orientation.

\begin{figure}
\centeredepsfbox{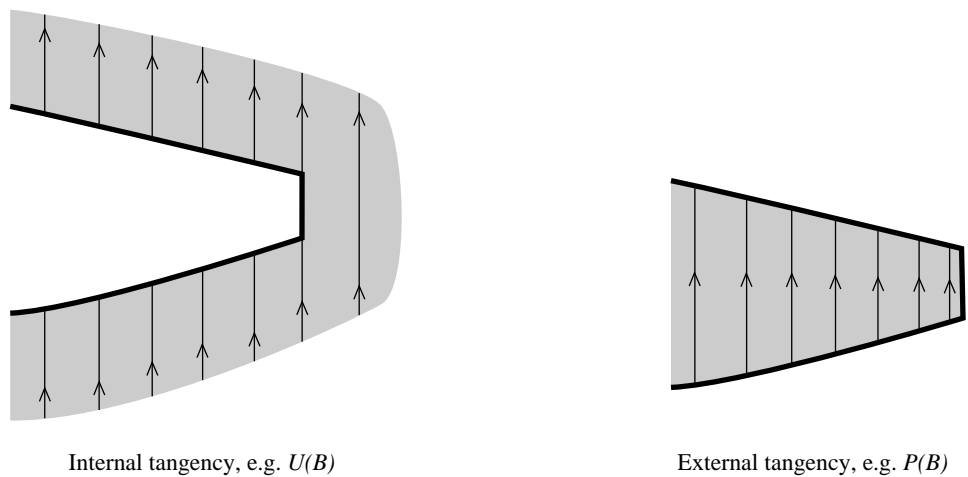}
\caption{Boundary tangencies of semiflows}
\label{BoundaryTan}
\end{figure}

Given a branched surface $B \subset M$, a {\em semiflow
neighborhood} $U(B)$ is a piecewise smooth regular neighborhood
equipped with a semiflow $\varphi$ satisfying the following conditions.
Each orbit of $\varphi$ is compact, transverse to $B$, and pierces $B$
in at least one point. There is a decomposition of $\bdy U(B)$ into
subsurfaces with disjoint interior
$$\bdy U(B) = \bdyv U(B) \union \bdyin U(B) \union \bdyout U(B)
$$
such that orbits of $\varphi$ point inward along $\bdyin
U(B)$, outward along $\bdyout U(B)$, and are {\em internally
tangent} along
$\bdyv U(B)$ as shown in figure \ref{BoundaryTan}. Each component of
$\bdyv U(B)$ is an annulus and the orbits of $\varphi$ define a fibration
of the annulus over the circle with interval fiber. If one forgets the
orientation and parameterization of $\varphi$, one obtains an $I$-bundle
neighborhood of $B$ as in \cite{Ga-Oe}.

We shall regard the {\em exterior} of $U(B)$ as a sutured manifold, as
follows.

\begin{figure}
\centeredepsfbox{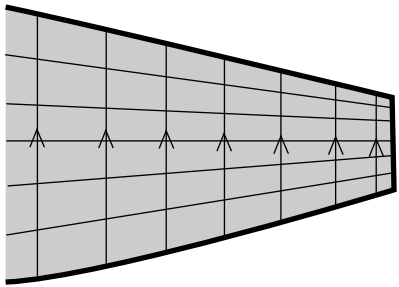}
\caption{A foliation and semiflow on a sutured manifold.}
\label{SuturedFoliation}
\end{figure}

A {\em sutured manifold} is an oriented 3-manifold $P$ whose boundary
is decomposed into surfaces with disjoint interior as $\bdy P = \R_- P
\union \R_+ P \union \Annuli P$, where $\R_- P$ and $\R_+ P$ are
disjoint, $\Annuli P$ is a collection of annuli called the
{\em sutures}, $\R_- P \intersect \Annuli P = \bdy(\R_- P) \subset \bdy
(\Annuli P)$, $\R_+ P \intersect \Annuli P = \bdy(\R_+ P) \subset \bdy
(\Annuli P)$, and each component of $\Annuli P$ has one boundary component
in $\R_- P$ and the other in $\R_+ P$. Let $\R_\pm P = \R_- P \union \R_+
P$. A sutured manifold $P$ is a {\em product} if there is a homeomorphism
$(P,\Annuli P) \homeo (S \cross I, \bdy S \cross I)$ for some compact
surface $S$. Given a semiflow on a sutured manifold $P$, we say that $P$
is an {\em isolating block} for the semiflow if orbits point outward
along $\R_+ P$, inward along $\R_- P$, and are externally tangent along
$\gamma P$ as shown in figure \ref{BoundaryTan}. A {\em foliation} of a
sutured manifold $P$ is required to be tangent to $\R_\pm P$ and
transverse to $\Annuli P$.  Near $\Annuli P$, a transverse foliation and
semiflow appear in cross section as in figure \ref{SuturedFoliation}.

If $B \subset M$ is a branched surface then $P(B) = \closure(M - U(B))$
has the structure of a sutured manifold where $R_+ P(B) = \bdyin U(B)$,
$R_- P(B) = \bdyout U(B)$, and $\Annuli P(B) = \bdyv U(B)$. Given a flow
$\varphi$ on $M$ transverse to $B$, there exists a regular neighborhood
$U(B)$ so that the restriction of $\varphi$ makes $U(B)$ into a semiflow
neighborhood, and so that $P(B)$ is an isolating block for $\varphi$. We
say that $U(B)$ is {\em adapted} to $\varphi$.

A {\em hierarchy} on $M$, in the sense of Gabai, is a sequence of
branched surfaces $B_0 \subset B_1 \subset \cdots \subset B_{N}$ such
that $B_0$ is a surface, $B_n - B_{n-1}$ is a union of sectors of
$B_n$ for each $1 \le n \le N$, and $P(B_{N})$ is a product. Let
$S(B_n) = B_n \intersect P(B_{n-1})$; this is a transversely oriented
surface properly embedded in $P(B_{n-1})$, and it is obtained from the
sectors $B_n - B_{n-1}$ by removing an open collar of each boundary
component. We may assume that each component $c$ of $\bdy S(B_n)$
satisfies the property that either $c$ is a core curve of some suture of
$P(B_{n-1})$ or $c$ is transverse to the sutures, meaning that each
component of $c \intersect \Annuli P(B_{n-1})$ is an arc connecting
opposite boundaries of some component of $\Annuli P(B_{n-1})$. This
property is immediate if each component of the boundary of $B_n - B_{n-1}$
either is in general position with respect to the branch locus of
$B_{n-1}$, or is contained in the branch locus. Otherwise, this property
can be achieved by first choosing $U(B_{n-1})$ so that the annuli of
$\bdyv U(B_{n-1})$ are extremely thin, and then doing a small isotopy of
$S(B_n)$. Note that $P(B_n)$ is obtained by doing a sutured manifold
decomposition of $P(B_{n-1})$ along $S(B_n)$, in the sense of \cite{Ga1}.

Let $\fol$ be a transversely oriented foliation in $M$, with a fixed
transversal flow $\varphi$. Given a saturated open set $W$, let $\widehat
W$ be its metric completion. The inclusion $\iota \from W \inject M$
extends to an immersion $\hat \iota \from \widehat W \rightarrow M$
carrying each component of $\partial \widehat W$ diffeomorphically onto a
leaf of $\fol$, but $\hat\iota$ may identify some of these boundary
components pairwise. In addition the structure of the foliation in
$\widehat W$ is as follows \cite{Di,Ca-Co}: $\widehat W = Q \union 
(\Sigma \times I)$, where $Q$ is a compact (possibly empty) sutured
manifold, each component of $\Sigma$ is a noncompact surface with compact
connected boundary, $Q$ is glued to $\Sigma \cross I$ by identifying
$\Annuli Q$ with $\bdy \Sigma \cross I$, and the restrictions of $\fol$
and $\varphi$ to $Q$ and to $\Sigma \cross I$ are well-behaved as
follows. The restrictions to $Q$ give a sutured manifold foliation and
semiflow; and the restriction of $\varphi$ to $\Sigma \cross I$ is a
product flow, that is, orbits have the form $x \cross I$. Note that
the restriction of $\fol$ to $\Sigma \cross I$ {\em need not} have
leaves of the form $\Sigma \cross t$. If a connected
open saturated set $W$ has $Q = \emptyset$, then $\widehat W$ (or $W$)
is called a {\em foliated  product}. An important fact is that for any
saturated open set $W$, at most finitely many components of $W$ are not
foliated products \cite{Di,Ca-Co}.  

The results above are also true for foliations of sutured manifolds. To
see this, first double along the sutures and then double along the
remaining boundary and apply the result to the final closed manifold.

\begin{proposition}{(hierarchy exists)}{} 
If $\Fol$ is a taut finite depth foliation on $M$ and $\varphi$ is a flow
transverse to $\Fol$, then there exists a hierarchy $B_0 \subset \cdots
\subset B_N$ transverse to $\varphi$ with neighborhoods $U(B_0) \subset
\cdots \subset U(B_N)$ adapted to $\varphi$, such that the sutured
manifolds $P(B_n)$ and the surfaces $S(B_n)$ are $\pi_1$-injective in
$M$, and each component of $B_0$ is isotopic to a compact leaf of $\Fol$.
\label{hierarchy}
\end{proposition}

The construction of finite depth foliations in theorem 5.1 of
\cite{Ga1} proceeds by first constructing a hierarchy and then
using it to produce the foliation. The point of this proposition is that all
finite depth foliations arise by this construction; and we obtain
additional information about flows. Since we cannot find a published
proof we provide one, though the ideas are well-known.

\begin{proof}{} Let $\Fol_n$ be the union of leaves of $\Fol$ of depth at
most $n$. Since $\fol$ is proper, $\fol_n$ is closed, hence it is a
sublamination of $\fol$, therefore an essential lamination in $M$. Let
$M_n = M - \fol_n$, and $\widehat M_n$ its metric completion. Then at
most finitely many components of $M_n$ are not foliated products. Now we
may alter $\fol$ without altering $\varphi$, collapsing foliated product
components of $M_n$; do this  inductively for each $n$. After the
alterations, only finitely many leaves are isolated, where a leaf $L$ at
depth $n$  is said to be isolated if $L$ is  isolated in $\fol_n$. In
addition if $L$ at depth $n$ is not isolated, then the component of
$\fol_n - \fol_{n-1}$  containing $L$ fibers over the circle with fiber
$L$. Next fatten up each of the finitely many isolated leaves into a
fibration over a closed interval, again altering $\fol$ without altering
$\varphi$. Each connected component of $\Fol_n - \Fol_{n-1}$ is now
either a fibration over a closed interval, or a fibration over the
circle; the latter are called ``circular'' components of $\Fol_n -
\Fol_{n-1}$.  Under these condition $\widehat M_n$ is a 1-1 immersed
submanifold of $M$ whose boundary components are leaves of $\Fol_n$. 

Now we construct by induction a hierarchy $B_0 \subset \cdots \subset
B_N$ and neighborhoods $U(B_0) \subset \cdots \subset U(B_N)$ adapted to
$\varphi$ so that if $0 \le n \le N$ then:

\begin{enumerate}
\item $U(B_n)$ contains every noncircular component of $\Fol_n -
\fol_{n-1}$, as well as an interval's worth of leaves in every circular
component.

\item $\Fol$ restricts to a foliation of the sutured manifold 
$P(B_n) = \overline{\bigl(M - U(B_n)\bigr)}.
$
It may be a product foliation
in some components of  $P(B_n)$.

\item The embedding $P(B_n) \inject P(B_{n-1})$ is $\pi_1$-injective.

\item $S(B_n) = B_n \intersect P(B_{n-1})$ is $\pi_1$-injective in
$P(B_{n-1})$.
\end{enumerate}

If $\Fol$ is not a fibration over $S^1$, let $U(B_0)$ to be the union of 
the compact leaves, otherwise let $U(B_0)$ be a closed interval of leaves.
Properties (1--2) are evident. We interpret (3) and (4) by setting
$P(B_{-1}) = M$, and these properties follow because the compact leaves
of $\Fol$ are incompressible surfaces in $M$. 

To continue the induction, given a component $Z$ of $P(B_n)$ on which
the restriction of $\Fol$ and $\varphi$ does not already give a product
structure, we describe the intersections with $Z$ of $P(B_{n+1})$,
$U(B_{n+1})$, and $B_{n+1}$ itself. Let $V(\bdy_\pm Z)$
be a neighborhood so that the restriction of $\phi$ has compact
orbits, defining a product structure with projection map $q \from
V(\bdy_\pm Z) \to \bdy_\pm Z$. 

Let $\fol'$ be the restriction of $\fol_{n+1}$ to $Z$. This lamination
has depth $\le 1$. Its compact leaves consist of $\bdy_\pm Z$ plus
foliated products $F \cross I$ with $\bdy F \cross I \subset \Annuli
Z$; all but finitely many of these foliated products have leaves
parallel to a component of $\bdy_\pm Z$. The noncompact leaves of
$\fol'$ fall into a finite collection of circular components and
foliated products, and each end $E$ of a noncompact leaf $L$
spirals into some component $F$ of $\bdy_\pm Z$ in the following
manner: there is a nonseparating, properly embedded, connected
1-manifold $c \subset F$ called a {\em juncture} for $E$, there is a
$\integers$-covering space
$\tilde F \to F$ to which $c$ lifts homeomorphically, there is a
subset $L_E \subset L$ representing $E$ and contained in $V(F)$ (the
component of $V(\bdy_\pm Z)$ containing $F$), and there is an embedding
$L_E \to\tilde F$, such that the maps $L_E \to \tilde F \to F$ and $L_E
\inject V(F) \to F$ are the same, and the image of $\bdy L_E$ is the
curve $c$. If $E, E'$ are ends of noncompact leaves spiralling into
the same component $F$ of $\bdy_\pm Z$, and if $E,E'$ are not in the
same end of a circular component of $\fol'$, we may assume that the
junctures for $E,E'$ are disjoint but isotopic curves in~$F$.

Let $\fol''$ be the union of the finitely many foliated products in
$\fol'$ which are not contained in $V(\bdy_\pm Z)$, plus a closed
interval of leaves in each circular component of
$\fol'$, plus $\bdy_\pm Z$. Since $Z - \fol''$ is an open subset of
$Z$ saturated by $\fol \restrictedto Z$, letting ${\overline {Z -
\fol''}}$ be the metric completion of $Z - \fol''$, there is a
decomposition ${\overline {Z - \fol''}} = Q \union (\Sigma \cross I)$
where $Q$ is a compact sutured manifold and $\Sigma \cross I$ is a
foliated product so that each noncompact component of $\Sigma$ has
compact, connected boundary. We can choose $Q$ big enough so that
$\Sigma \cross I$ is contained in $V(\bdy_\pm Z)$, and so that the
boundary of each component of $\Sigma$ projects to some juncture, thus
we can associate a juncture to each component of $\Sigma$. Moreover, if
$\Sigma_1,\ldots,\Sigma_n$ are components of $\Sigma$ with isotopic
junctures $c_1,\ldots, c_n$ lying in a component $F$ of $\bdy_\pm
Z$, then we can choose $Q$ so that the picture in figure
\ref{junctures} holds: the curves $c_1,\ldots,c_n$ are ordered so
that $c_i \union c_{i+1}$ bounds an annulus disjoint from the other
curves; and the ``first hitting'' map $\Sigma_{i+1} \to \Sigma_{i}$ is
continuous for $i=1,\ldots,n-1$, where this map is defined by starting
from a point in $\Sigma_{i+1}$, going along an orbit of $\phi$ towards
$F$, and stopping at the first point of $\Sigma_{i}$.

\begin{figure}
\centeredepsfbox{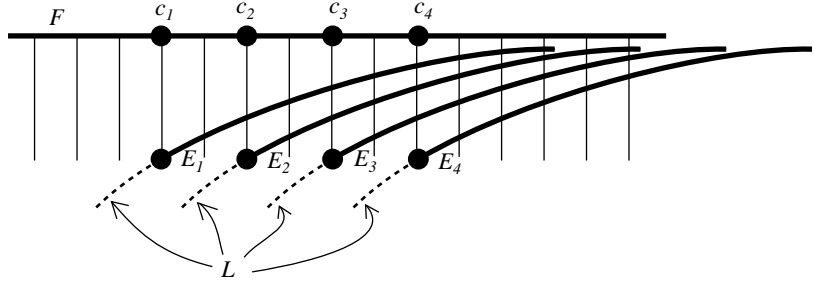}
\caption{Isotopic junctures}
\label{junctures}
\end{figure}

Now define $P(B_{n+1}) \intersect Z = Q$. Notice that $\fol''$ is an
essential lamination of $Z$, hence $Z - \fol''$ is
$\pi_1$-injective in $Z$. Since each noncompact component of $\Sigma$
has compact, connected boundary, it follows that $Q$ is
$\pi_1$-injective in $Z - \fol''$, and property (3) for
$P(B_{n+1})$ follows.  Property (2) is easily checked, as is property
(1) for $U(B_{n+1}) = {\overline {(M - P(B_{n+1}))}}$.

To define the sectors of $B_{n+1}$ intersecting $Z$, note that the
foliation and flow define a product structure
$\fol'' - \bdy_\pm Z = G \cross [0,1]$ for some surface $G$, not
necessarily connected. Each end of $G \cross 0$ or $G \cross 1$ is
eventually contained in $\Sigma \cross I$. Let $L$ be the set of $x \in
G$ such that $x \cross [0,1]$ is disjoint from $\Sigma \cross I$.
Then $L$ is a subsurface of $G$, and we may regard $L$ as being
embedded in $\fol''$ transverse to the flow. Homotop $\bdy L$ to a map
into $B_n$ as described below, and extend the homotopy over $L$. The
image of $L$ after the homotopy is the union of sectors of $B{n+1} -
B_n$ intersecting $Z$. To describe the homotopy of $\bdy L$, each
component $\beta$ of $\bdy L - \Annuli Z$ lies in
$V(\bdy_\pm Z)$, and $\beta$ projects via $\phi$ to a juncture;
homotop $\bdy L$ along orbits of $\phi$ until $\beta$ lies on that
juncture. Extend the homotopy of $\bdy L - \Annuli Z$ over all of
$\bdy L$, so that each arc component of $\bdy L \intersect \Annuli Z$
maps to an arc in $B_n$ crossing the branch locus in a single
transverse intersection point, and each circle component of $\bdy L
\intersect \Annuli Z$ maps to a circle in the branch locus of $B_n$.
Because of the properties of junctures given in figure
\ref{junctures}, this homotopy may be carried out so that $L - \bdy
L$ is embedded in $M - B_n$, and hence the image of $L$ after the
homotopy consists of sectors of $B_{n+1} - B_n$ whose
boundary lies in $B_n$, as required for a hierarchy. 

Property (4) for $S(B_{n+1})$ follows because each component of
$S(B_{n+1})$ corresponds to a component of $L$ as above, but
components of $L$ are $\pi_1$-injective in in $P(B_n)$.

The only nonobvious points remaining in the proof of the lemma are the
statements on $\pi_1$-injectivity. Using (3) it follows by induction that
$P(B_n)$ is $\pi_1$-injective in $M$. Using (4), and the fact that
$P(B_{n-1})$ is $\pi_1$-injective in $M$, it follows that $S(B_n)$ is
$\pi_1$-injective in $M$.
\end{proof}

We also need some geometric information about hierarchies. Consider a
closed, oriented 3-manifold $M$ and a nonseparating incompressible
surface $F$. We say $F$ is a {\em fiber} if $F$ is a leaf of a
fibration of $M$ over $S^1$. If $M$ is hyperbolic we say that $F$ is
{\em quasi-fuchsian} if the representation $\pi_1 F \to \pi_1 M \to
\Isom(\hyperbolic^3)$ is quasi-fuchsian. A deep fact due to Marden
\cite{Ma}, Thurston \cite{Th1} and Bonahon \cite{Bo}, is that when
$M$ is hyperbolic then $F$ is either a fiber or quasi-fuchsian (if
$F$ is separating there is another option, namely that $F$ is a
``virtual fiber'').

Recall from \cite{Gr,Gh-Ha} that a locally compact path metric space
$X$ is {\em negatively curved in the large} if it satifies the ``thin
triangles'' condition: there exists $\delta>0$ such that for every
geodesic triangle with sides $A_1,A_2, A_3$  we have $A_1
\subset U_\delta(A_2 \union A_3)$. In this case there is a
canonically defined ideal boundary at infinity denoted
$\bdy_\infinity X$, and there is a canonical compactification $X
\union \bdy_\infinity X$. If $X$ is negatively curved in the large,
so is any space quasi-isometric to $X$. If $X,Y$ are negatively
curved in the large and $f \from X \to Y$ is a quasi-isometry then
$f$ extends canonically to a continuous map $\hat f \from X \cup
\bdy_{\infty} X \rightarrow Y \cup \bdy_{\infty} Y$ that
restricts to an embedding $\bdy_\infinity X \inject \bdy_\infinity Y$.
For example, hyperbolic space $\hyperbolic^n$ is negatively curved in the
large, and its ideal boundary in the sense of Gromov is the same as the
sphere at infinity $S^{n-1}_\infinity$. Also, any closed convex subset $X
\subset \hyperbolic^n$ is negatively curved in the large, and
$\bdy_\infinity X$ coincides with the set of limit points of $X$ in
$S^{n-1}_\infinity$. 

Given a group $G$ with a finite generating set $A$, the {\em Cayley
graph} $C$ is the graph with a vertex for each $g \in G$, and an edge from
$g$ to $ga$ for each $a \in A$. By taking each edge to be a path of
length 1 we make the Cayley graph a path metric space. For any two finite
generating sets, the Cayley graphs are quasi-isometric. The group $G$ is
{\em word hyperbolic} if some (and hence any) Cayley graph of $G$ is
negatively curved in the large. 
Then we let $\partial _{\infty} G = \partial _{\infty} C$.
Consider a compact Riemannian
manifold $M$ with $G = \pi_1(M,x)$, choose a lift $x_0$ of $x$ in the
universal cover $\widetilde M$ and choose loops in $M$ representing a
finite generating set for $\pi_1 (M)$. Lifting to $\widetilde M$ we have a
map of the Cayley graph of $G$ to $\widetilde M$, which gives a 
quasi-isometry between $G$ and $\widetilde M$ \cite{Mi}. Therefore $G
= \pi_1 M$ is word hyperbolic if and only if $\widetilde M$ is
negatively curved in the large, in which case there is a canonical
identification
$\bdy_\infinity G \homeo \bdy_\infinity \widetilde M$.

\begin{proposition}{}{} Let $M$ be a closed, orientable hyperbolic
3-manifold, $B_0 \subset \cdots \subset B_N$ a hierarchy in $M$ such
that $P(B_n)$ and $S(B_n) = B_n \intersect P(B_{n-1})$ are
$\pi_1$-injective for all $n$. Suppose that some component $F$ of
$B_0$ is a quasi-fuchsian surface in $M$. Then:

\smallskip
\noindent
(1) Each group $\pi_1(P(B_n))$, $\pi_1(S(B_n))$ is word hyperbolic.

\noindent 
(2) Each embedding $\cover S(B_n) \inject \cover P(B_{n-1})$ is a
quasi-isometry which takes $\bdy \cover S(B_n)$ into
$\bdy \cover P(B_n)$. Therefore it extends to an embedding $\cover
S(B_n) \union \bdy_{\infty} \cover S(B_n) \inject \cover P(B_{n-1})
\union \bdy_\infinity \cover P(B_{n-1})$.

\noindent
(3) Each inclusion $\cover P(B_{n}) \inject \cover P(B_{n-1})$ is a
quasi-isometry in the path metrics induced from $\mi$.
Therefore it extends to an embedding
$\cover P(B_n) \cup \partial _{\infty} \cover P(B_n) \inject
\cover P(B_{n-1}) \cup \partial _{\infty} \cover P(B_{n-1})$.

\label{DiscsAndBalls}
\end{proposition}

\noindent
Remark: The group $\pi_1 (S(B_n))$ may be cyclic, when $S(B_n)$
is an annulus; or trivial, when $S(B_n)$ is a disc.

\begin{proof}{}  
Let $\rho \from \pi_1(M) \rightarrow PSL_2({\bf C})$ be the holonomy
representation of the hyperbolic structure. Consider the subgroup $G =
\rho(\pi_1(M-F))$ with limit set $\Lambda G \subset S^2_\infinity$, let
$\Hull (\Lambda G) \subset \hyperbolic^3 \union S^2_\infinity$ be its convex
hull, let $\Hull_\epsilon (\Lambda G) \subset \hyperbolic^3 \union
S^2_\infinity$ be the union of $\Hull (\Lambda G)$ with its
$\epsilon$-neighborhood in $\hyperbolic^3$, and let $X_G =
\Hull_\epsilon (\Lambda G) \intersect \hyperbolic^3$. Since
$F$ is quasi-Fuchsian, the group $G$ is convex cocompact, that is
$G$ acts cocompactly on $X_G$; equivalently, $G$ is geometrically finite
\cite{Ma,Th1}. It follows that any finitely generated subgroup $H
\subgroup G$ is convex cocompact \cite{Th1,Cn}. Since $X_H$ is negatively
curved in the large, and the action of $H$ on $X_H$ is properly
discontinuous and cocompact, it follows that $H$ is word hyperbolic,
proving (1).

Suppose $P \subset M-F$ is a compact, $\pi_1$-injective submanifold and $H
= \pi_1 (P)$. We show that the natural inclusion $j \from \cover P \to
\hyperbolic^3$ is a quasi-isometric embedding, where $\widetilde P$ 
has the path metric induced from $\mi$.

The natural inclusion
$\iota \from X_H \inject \hyperbolic^3$ is a quasi-isometric embedding.
Let $\Gamma_H$ be the Cayley graph of $H$, and consider the two
quasi-isometries $\alpha \from\Gamma_H \to X_H$, $\beta \from \Gamma_H
\to
\cover P$ given by \cite{Mi}. The maps $\iota\composed \alpha, j
\composed \beta \from \Gamma_H \to\hyperbolic^3$ differ by a bounded
amount, therefore $j \composed \beta$ is a quasi-isometric embedding.
There is a quasi-isometry $\bar\beta \from \cover P \to \Gamma_H$
which is an inverse of $\beta$ in the quasi-isometric category, so in
particular
$\beta \composed \bar\beta \from \cover P \to \cover P$ differs from
the identity by a bounded amount.  Thus we have a quasi-isometric
embedding
$j \composed \beta \composed \bar\beta \from \cover P \to
\hyperbolic^3$ which differs from $j$ by a bounded amount, so $j$ is a
quasi-isometric embedding.  Hence each $\widetilde P(B_n) \inject
\hhh$ extends to a  continuous embedding $\cover P(B_n) \cup \partial
_{\infty} 
\cover P(B_n)  \inject \hhh \cup \si$.

To prove (3), the quasi-isometric embedding $\cover P(B_n) \inject
\hyperbolic^3$ factors through the natural embedding $\cover P(B_n)
\inject \cover P(B_{n-1})$, so $\cover P(B_n) \inject \cover
P(B_{n-1})$ is a quasi-isometric embedding. The proof of (2) is
similar.
\end{proof}

\smallskip
Remark: Although not logically necessary for our
results, it is helpful to keep in mind the following additional facts:

\smallskip
\noindent
(4) Each compactified universal cover $\widetilde P(B_n) \cup
\bdy_{\infty} \widetilde P(B_n)$ is a 3-ball.

\noindent
(5) Each $\widetilde S(B_n) \union \bdy_{\infty} \widetilde
S(B_n)$ is a 2-disc properly embedded in the above 3-ball. 
\smallskip

To prove (4), let $P = P(B_n)$, let $H = \pi_1(P)$, and consider
the action of $H$ on $\hyperbolic^3$. In the manifold $\hyperbolic^3 /
H$, both of the compact manifolds $X_H / H$ and $P$ embed as deformation
retracts. By \cite{Mc-MS} it follows that there is a homeomorphism $X_H /
H \homeo P$ in the correct homotopy class. This homeomorphism lifts to an
$H$-equivariant quasi-isometric homeomorphism $X_H \homeo \cover P$, which
extends to a homeomorphism $\Hull_\epsilon (\Lambda H) \homeo \compact P$,
and $\Hull_\epsilon (\Lambda H)$ is obviously a 3-ball.

The proof of (5) when $S(B_n)$ is a disc or annulus is easy. Otherwise
there is a hyperbolic metric with geodesic boundary on $S(B_n)$, and
$\cover S(B_n) \union \bdy_\infinity S(B_n)$ is therefore a 2-disc,
which by (2) is properly embedded in the 3-ball $\compact P$.

\newcommand{\flowtime}{\func{\tau}}
\newcommand{\Fr}{\func{Fr}}
\newcommand{\Bc}{\func{Bc}}
\newcommand{\ClFr}{\sigma^f}
\newcommand{\ClBc}{\sigma^b}
\newcommand{\Axis}{\func{Axis}}

\section{Inductive proof of theorem A}
\label{qg}

To set up the induction, apply proposition \ref{hierarchy} to obtain a
hierarchy $B_0\subset\cdots\subset B_N$. We use the following
notation. If $0 \le n \le N$ let $P_n = P(B_{N-n})$, and let $P_{N+1} =
M$, so $P_0 \subset P_1 \subset \cdots \subset P_N \subset P_{N+1}$;
the indexing is reversed to facilitate the induction proof. If $1 \le
n \le N+1$ let $S_n = S(B_{N-n+1}) = B_{N-n+1} \intersect P_n$, so
$S_n$ is properly embedded in $P_n$, and $P_{n-1}$ is obtained from
$P_n$ by removing a regular neighborhood of $S_n$. Since some (and
hence all) compact leaves of $\Fol$ are quasi-fuchsian, the same is
true for components of $B_0$, hence proposition \ref{DiscsAndBalls}
applies.

Given a flow line $x \cdot \reals$ of $\Phi$, define the {\em codepth}
to be the minimal integer $n$ such that $x \cdot \reals \subset P_n$,
and define $\Omega_n$ to be the union of all flow lines of codepth at
most $n$. Note that $\Omega_n$ is the set of all flow lines contained
in $P_n$, and $\Omega_n$ is closed. As special cases, $\Omega_0 =
\emptyset$ since $P_0$ is a product, and $\Omega_{N+1} = M$. 

The program for the proof of theorem A
is to assume that orbits in $\Omega_{n-1}$ 
are uniformly quasigeodesic and then show that orbits in $\Omega_n$
are also uniformly quasigeodesic (with a bigger quasigeodesic
constant). Since $\Omega_0 = \emptyset$ and $\Omega_{N+1} = M$,
induction will show that all orbits are uniformly quasigeodesic.
This will complete the proof of theorem A.

For notational convenience, throughout this section we write $P =
P_n$, $\Omega = \Omega_n$, $S = S_n$, $P' = P_{n-1}$, and $\Omega' =
\Omega_{n-1}$. Let $\pi \from \cover P \to P$ be the universal
covering, and let $\cover \Omega = \pi^\inverse(\Omega)$. 

Fix a connected lift $\cover P' \inject \cover P$, and let $\cover
\Omega' = \pi^\inverse(\Omega') \intersect \cover P'$. The induction
hypothesis says that orbits in $\cover \Omega'$ are uniform
quasigeodesics in $\cover P'$. Since $\cover P' \inject \cover P$ is a
quasi-isometry, then orbits in $\cover \Omega'$ are uniformly
quasigeodesic in $\cover P$; using the action of $\pi_1 P$ by
isometries on $\cover P$, the same is true for orbits in
$\pi^\inverse(\Omega')$.
Recall that the hypothesis in theorem $A$ is
that $\wwp$ has Hausdorff orbit space.
This will be used in verifying conditions (b) and (c) of
theorem \ref{upgrade}.

We will need the following well known simple result \cite{Co}:

\begin{lemma}{}{} Let $W$ be a compact metric space with a nonsingular
semiflow $\varphi$ parameterized by arc length. Let $\Omega$ be the set
of points $x$ for which $\varphi_t(x)$ is defined for all $t \in
\rrrr$. Then given any $\delta > 0$, there is $a > 0$ so that any orbit
$\gamma$ of $\varphi$ is in the $\delta$-neighborhood of $\Omega$
except perhaps for an initial segment of length $< a$ and another
final segment of length $< a$.
\end{lemma}

We will also use the following localization property of quasigeodesics.

\begin{proposition}{}{\cite{Gr,Gh-Ha}} Let $M$ be a compact manifold 
with negatively curved
$\pi_1(M)$. Then for any $K > 0$ there is $L(K)  > 0$ (usually $L(K) >>
K$) satisfying: if $\gamma$ is an embedded curve so that any subarc of
length $\leq L(K)$ is a $K$ quasigeodesic then $\gamma$ is a $2
K$-quasigeodesic.
\label{local}
\end{proposition}

The following essential fact which is a direct
consequence of proposition \ref{DiscsAndBalls} will
be used explicitly or implicitly throughout this section:
if $\gamma$ is a curve contained in $\widetilde P'$, then
$\gamma$ is a quasigeodesic in $\widetilde P'$ if and only if
$\gamma$ is a quasigeodesic in $\widetilde P$ and also
 if and only if
$\gamma$ is a quasigeodesic in $\hhh$.
The quasigeodesic constants may differ.
We caution the reader that some arguments are done in $\hhh$
while others are done in $\widetilde P$. The context makes it clear.

For flow segments disjoint from $S$, the following lemma establishes
the quasigeodesic property directly.

\begin{lemma}{}{}
There is $K > 0$ so that all flow segments, half orbits or full orbits
of $\wwp$ contained in $\cover P'$ are $K$-quasigeodesics of $\cover
P$; translating by the action of $\pi_1(P)$, the same is true in any
lift of $P'$.
\label{uniform}
\end{lemma}

\begin{proof}{}
The full orbits staying in $\cover P'$ are precisely the orbits in
$\cover\Omega'$, which are $k$-quasigeodesic for some uniform $k$. Fix
$\delta_0 > 0$, with $2 k \delta_0 < 1$. Let $L = L(k+1)$ be given by
proposition \ref{local}. Choose $\delta_1 > 0$ so that if $x, y \in
\cover P$ and $d(x,y) < \delta_1$ then $d(x \cdot t, y \cdot t) <
\delta_0$ for any $|t| < L$. Choose $a > 0$ so that any orbit in $P'$
is in the $\delta_1$-neighborhood of $\Omega'$, except perhaps for
initial and final segments of length $< a$.

Let now $x, y \in \cover P'$ with $y = x \cdot t$. Choose a
subsegment $\gamma = [z,w] \subset [x,y]$, with $z = x \cdot t_1$, $w
= y \cdot t_2$, $0 \le t_1 < a$ and $-a < t_2 \le 0$, so
that $\gamma \subset U_{\delta_1} (\cover \Omega')$. Let $\alpha =
[z_0,w_0]$ be a subarc of $\gamma$ with $0 < t_3 = \flowtime(z_0,w_0)
\le L$. By the choice of $\delta_1$ it follows that  there is $z_1 \in
\cover \Omega'$, with $d(z_0 \cdot t, z_1 \cdot t) < \delta_0$, for
any $|t| < L$. Since $z_1 \cdot \reals$ is a $k$-quasigeodesic then
\begin{eqnarray*}
\flowtime(z_0,w_0) = t_3 = \flowtime(z_1,z_1 \cdot t_3) 
\leq k \, d(z_1,z_1 \cdot t_3) + k \\
\leq k d(z_0,w_0) + 2 k \delta_0 +
k < k d(z_0,w_0) + (k + 1).
\end{eqnarray*}
Clearly this also works for any subsegment of $\alpha$, hence $\alpha$
is a $(k+1)$-quasigeodesic. By the previous proposition, one concludes
that $\gamma$ is a $(2k+2)$-quasigeodesic. This implies that 
\begin{eqnarray*}
\flowtime(x,y) & < 2a + \flowtime(z,w) \le (2k+2) d(z,w) + (2k+2+2a) \\
     & \le (2k+2) \bigl( d(x,y) + 2a \bigr) + (2k+2+2a) \\
     & = (2k+2) d(x,y) + \bigl( 2k+2+2a + (2k+2)2a \bigr ) 
\end{eqnarray*}
Hence any piece of orbit of $\wwp$ contained in $\cover P'$ is
a $(4ka + 2k + 6a + 2)$-quasigeodesic of $\widetilde P'$.
Since $\widetilde P'$ is quasi-isometrically embedded in
$\widetilde P'$, there is $K > 0$ so that any
piece of orbit of $\wwp$ contained in $\widetilde P$ is a 
$K$-quasigeodesic 
of $\widetilde P$.
\end{proof}

Now we prepare the ground for applying theorem \ref{upgrade} to
show that full orbits in $\cover\Omega$ are uniformly quasigeodesic. We
must study how orbits in $\cover\Omega - \cover\Omega'$ cross lifts of
$S$ in $\cover P$, so we embark on a study of these lifts.

Let ${\mathcal E}$ be the collection of lifts of $S$ to $\cover P$. Any
element $E \in {\mathcal E}$ is transversely oriented and separates
$\cover P$. The front of $E$ is the component $\Fr(E)$ of $\cover P -
E$ on the positive side of $E$, and the closure of $\Fr(E)$ in
$\widetilde P \cup \partial _{\infty} \widetilde P$,
is denoted $\ClFr(E)$. The back $\Bc(E)$ and its closure
$\ClBc(E)$ are similarly defined. Define a strict partial order on
$\mathcal E$ where $E < E'$ if $\Bc(E) \intersect \Fr(E') = \emptyset$ and
$E \not = E'$, hence in particular $E' \subset \Fr(E)$.

Notice that since $S$ is compact, a bounded subset of $\cover P$
intersects only finitely many $E \in {\mathcal E}$. The following
lemmas strengthen this fact, by showing that sequences in
$\mathcal E$ are limited in how they may accumulate in the ball 
$\widetilde P \cup \partial _{\infty} \widetilde P$;
these lemmas will be useful in analyzing flow lines that cross $S$ many
times.
If $\rho: \pi_1(M) \rightarrow PSL_2({\bf C})$ is the holonomy
representation, then $\rho(\pi_1(P))$ is a Kleinian group
which is convex cocompact.

\begin{lemma}{}{}
There is $J_0 > 0$ such that if $E_1, \ldots, E_{J_0}$ are distinct
elements of $\mathcal E$, then $\bdy_\infinity(E_1) \intersect\cdots
\intersect \bdy_\infinity(E_{J_0}) = \emptyset$. If moreover $E_1 <
\cdots < E_{J_0}$ then $\sigma^b E_1 \intersect \sigma^f E_{J_0} =
\emptyset$.
\label{LimitSetsDisjoint}
\end{lemma}

\begin{proof}{}
First we need the fact that if $H_1, \ldots, H_n$ are geometrically
finite Kleinian groups that generate a discrete group, then $\Lambda
H_1 \intersect \cdots \intersect \Lambda H_n = \Lambda(H_1 \intersect
\cdots \intersect H_n)$. When $n=2$ this is proved in \cite{Su},
and the statement for finitely
many subgroups follows by induction.

To prove the first statement, suppose that $\bdy_\infinity(E_1)
\intersect\cdots \intersect \bdy_\infinity(E_J) \ne \emptyset$. Let
$H_i\subset \pi_1 P$ be the stabilizer subgroup of $E_i$, so
$\Lambda(H_i) = \bdy_\infinity(E_i)$, and by the above argument using
Susskind's theorem it follows that $H_1 \intersect \cdots \intersect
H_J \ne \emptyset$. Let $f$ be in the intersection, and let $\Axis_f$
be an axis for $f$ in $\hyperbolic^3$. By conjugation, we may assume
that $\Axis_f$ intersects a fixed fundamental domain $D$ for $\pi_1
M$. 

Since $E$ is quasi-isometrically embedded in $\hhh$, it
is $R$-quasiconvex, that is, for any 
$x, y \in  E \cup \partial _{\infty} E$,
the geodesic in $E$ connecting them is at most $R$ distant
from the hyperbolic geodesic connecting them \cite{Gr,Gh-Ha}.
The $R$ is independent of the lift $E$ of $S$.
Thus each $E_i$ intersects $U_R(D)$, the open $R$-neighborhood
of $D$. This shows that $J \le J_0$ where $J_0$ is an upper bound for
the number of distinct $E \in {\mathcal E}$ intersecting the bounded set
$U_R(D)$.

To prove the second statement, suppose that $\sigma^b E_1 \intersect
\sigma^f E_{J_0} \ne \emptyset$. 
It follows that $\bdy_\infinity E_1
\intersect \bdy_\infinity E_{J_0} \ne \emptyset$; let $\xi$ be a point
in this intersection. 
If for some $1 \leq i \leq n$, $\xi \not \in \bdy_{\infty} E_i$,
let $V$ be a neighborhood of $\xi$ in 
$\widetilde P \cap \partial _{\infty} \widetilde P$ with
$E_i \cap (V \cap \widetilde P) = \emptyset$.
However since $\xi \in \bdy_{\infty} E_0 \cap \bdy_{\infty} E_{J_0}$,
there are $x \in E_0 \cap V$ and $y \in E_{J_0} \cap V$.
Hence $x$ can be connected to $y$ in $V \cap \widetilde P$,
contradicting the fact that $E_i$ separates $E_0$ from $E_{J_0}$.
But then we have proved $\bdy_\infinity
E_1 \intersect \cdots\intersect \bdy_\infinity E_{J_0} \ne \emptyset$,
contradicting the first statement of the lemma.
\end{proof}

The $J_0$ given by the previous lemma is fixed from now on. The next
lemma gives an even stronger accumulation property.


\begin{lemma}{}{}
Given an infinite sequence $E_1, E_2, \ldots \in {\mathcal E}$ such that
$E_i\ne E_j$ if $i \ne j$, suppose there exists $E_0 \in {\mathcal E}$ such
that $E_0 < E_i$ for all $i \ge 1$. Then there is a subsequence
$E_{i_n}$ such that $\ClFr(E_{i_n})$ converges to a single point $\xi
\in\bdy_\infinity \cover P$, in the Hausdorff topology on closed
subsets of $\widetilde P \cup \partial _{\infty} \widetilde P$.
\label{BoundedBelow}
\end{lemma}

\begin{proof}{} Choose covering translations $f_i \in \pi_1 P$ such
that $f_i(E_0) = E_i$. Then $f_i \ne f_j$ when $i \ne j$, so by the
convergence group property for Kleinian groups applied to
$\rho(\pi_1(M))$ (see pg. 22 of Maskit's book
\cite{Ma}), we may pass to a subsequence so that there is a pair of
points $\xi_-, \xi_+ \in \bdy_\infinity P$ called a {\em source} and
{\em sink} for the sequence of functions $f_i$, meaning that $f_i
\restrictedto \bigl( \widetilde P \union \bdy_\infinity \widetilde P
\bigr) - \xi_-$ converges uniformly on compact sets to the constant
map with value $\xi_+$, and $f_i^\inverse \restrictedto \bigl(
\widetilde P \union \bdy_\infinity \widetilde P
\bigr) - \xi_+$ converges similarly to $\xi_-$. Since
$f_i^\inverse(\ClBc(E_i)) = \ClBc(E_0)$ for all $i$, it follows
that $f_i^\inverse(\ClBc(E_0)) \subset \ClBc(E_0)$, and hence $\xi_-
\in \ClBc(E_0)$. 

If $\xi_- \in \sigma^b E_0 - \bdy_\infinity E_0$ then we are done,
because then $\xi_- \notin \ClFr(E_0)$, so the sequence of maps $f_i
\restrictedto \ClFr(E_0)$ converges uniformly to the constant map with
value $\xi_+$, and hence the images $\ClFr(E_i)$ of these maps converge
in the Hausdorff topology to $\xi_+$.

Suppose on the other hand that $\xi_- \in \bdy_\infinity E_0$. By the
previous lemma, for all but finitely many $E_i$ we have $\xi_- \notin
\bdy_\infinity E_i$; in particular this is true for some $i = i_0$.
We also have $\xi_- \notin \ClFr(E_{i_0})$, because $\xi_- \in
\sigma^b(E_0)$. Note that $\xi_-, \xi_+$ is a source, sink pair for the
sequence $f_i \composed f_{i_0}^\inverse$, so the maps $f_i \composed
f_{i_0}^\inverse \restrictedto \ClFr(E_{i_0})$ converge uniformly to
the constant map with value $\xi_+$, and hence their images
$\ClFr(E_i)$ converge in the Hausdorff topology to $\xi_+$.
\end{proof}

Now we show that all flow lines in $\cover \Omega$ extend continuously
to $\bdy_\infinity \widetilde P$, verifying condition (a) of theorem
\ref{upgrade}.

\begin{proposition}{(extension of orbits)}{}
If $x \in \cover \Omega$, then the following limits exist:
$$\lim_{t \rightarrow +\infty} x \cdot t = \eta_+(x) \in
\bdy_\infinity \cover P \quad {\text and } \quad
\lim_{t \rightarrow -\infty} x \cdot t = \eta_-(x) \in
\bdy_\infinity \cover P
$$
\label{ext}
\end{proposition}

\begin{proof}{} We only consider forward limits. There are two cases:

\smallskip
\noindent
{\underline {Case 1}} --- $x \cdot \Rplus$  eventually stops
intersecting $\pi^{-1}(S)$ (this includes $x \in \cover \Omega'$).

Let $y$ be the last intersection of $x \cdot \RPlus$ with
$\pi^{-1}(S)$ if any exists, otherwise let $y = x$. Then $\pi(y) \cdot
\RPlus$ is contained in $P'$ hence $y \cdot \RPlus$ is a
$K$-quasigeodesic in $\cover P$ by lemma \ref{uniform}. Therefore it
has a unique limit point.

\smallskip
\noindent
{\underline {Case 2}} --- 
$x \cdot \RPlus$ keeps intersecting $\pi^{-1}(S)$.

 Let $E_0 < E_1 < \cdots \in {\mathcal E}$ be the elements that $x
\cdot \RPlus$ intersects. Thus, all accumulation points of $x
\cdot \RPlus$ are contained in
$\ClFr(E_1) \intersect \ClFr(E_2) \intersect \cdots $. Since this is
a nested intersection, it is the same as the Hausdorff limit of the
sets $\sigma^f(E_1), \sigma^f(E_2), \cdots$, which by
lemma \ref{BoundedBelow} is a single point.
\end{proof}

Next we verify condition (b) of theorem \ref{upgrade}.
We use the following facts: if $\alpha$ and $\beta$ are
two quasigeodesic rays in $\hhh$ with the same ideal point,
then there is $R > 0$ (depending on $\alpha$ and $\beta$)
so that $\alpha$ is in the $R$ neighborhood of $\beta$ and
vice versa. For a pair of bi-infinite quasigeodesics with 
the same ideal points there is also a bound, which
depends only on the quasigeodesic constant.
We also need the following very useful result:

\begin{lemma}{}{}
Suppose that $\gamma_i = x_i \cdot [0,t_i)$ and $\gamma = x \cdot
\Rplus$ are $k_1$-quasigeodesic flow segments or rays, where $x_i \to
x$ in $\cover P$ and $t_i \to +\infinity$ in $[0,+\infinity]$; we allow
$t_i = +\infinity$. Then $x_i \cdot t_i \to\eta_+(x)$ in 
$\widetilde P \cup \partial _{\infty} \widetilde P$,
where $x_i \cdot t_i$ means $\eta_+(x_i)$ if $t_i = +\infinity$.
\label{conv}
\end{lemma}

\begin{proof}{}
Otherwise, passing to a subsequence we have $x_i \cdot t_i \to w \in
\widetilde P \cup \partial _{\infty} \widetilde P$
with $w \ne \eta_+(x)$. Choose $R > 0$ so that any
$k_1$-quasigeodesic is at most $R$ distant from a corresponding minimal
geodesic. 

Choose disjoint small neighborhoods $W_0$ of $\eta_+(x)$ in 
$\compact P$,
$W_1$ of $w$ in $\compact P$ and $W_2$ of $x$ in $\cover P$ so that
for any minimal geodesic segment or ray $\beta$ starting in $W_2$ and
with endpoint in $W_1$ then

$$U_R(\beta) \cap (W_0 \cap \cover P) = \emptyset.$$ 

\noindent
For $i$ big enough $x_i \in W_2$ and $x_i \cdot t_i \in W_1$, so
since $\gamma_i$ is at most $R$ distant from the corresponding minimal
geodesic, then $\gamma_i \cap (W_0 \cap \cover P) = \emptyset$.
But since $x_i \rightarrow x$ then by continuity of the flow 
$\gamma_i \cap (W_0 \cap \cover P) \ne \emptyset$ for $i$ big
enough, contradiction.
\end{proof}

\begin{proposition}{(distinct limit points)}{} For any $x \in 
\cover \Omega$, $\eta_+(x) \not = \eta_-(x)$.
\label{distinct}
\end{proposition}

\begin{proof}{}
There are two cases:

\smallskip
\noindent
{\underline {Case 1}} --- $x \cdot \reals$ intersects $\pi^{-1}(S)$
only finitely many times. 

Then there are $t_0 < t_1 \in \rrrr$, with
$\gamma_0 = x \cdot (-\infty,t_0]$ and $\gamma_1 = x \cdot
[t_1,+\infty)$, so that $\gamma_0 \cap \pi^{-1}(S) = \gamma_1 \cap
\pi^{-1}(S) = \emptyset$. By lemma \ref{local}, $\gamma_0$ and
$\gamma_1$ are $K$-quasigeodesics in $\hhh$. Assume they have the same
ideal point $p$. Then there are sequences $v_i \in \gamma_0$, $u_i \in
\gamma_1$, such that $u_i \to p$, $v_i \to p$, and $d(u_i,v_i)$ is
bounded. But $\tau(u_i,v_i) \to \infinity$.

Let $y = x \cdot t_0$ and $z = x \cdot t_1$.
Up to taking a subsequence assume that $\pi(u_i)$ 
converges in $M$. Then there are covering translations $h_i$
so that $h_i(u_i) \rightarrow u$. Since $u_i \in \widetilde P'$,
then $\pi(u) \in P'$.
We can also assume that $h_i(v_i) \rightarrow v$ and hence
$\pi(v) \in P'$.

Since $[z,u_i]$ is contained in $\widetilde P'$,
it is a $K$-quasigeodesic by lemma \ref{uniform}.
Since $\tau(z,u_i) \rightarrow +\infty$,
it then follows that $d(z,u_i) \rightarrow +\infty$,
so $d(h_i(z), h_i(u_i)) \rightarrow +\infty$.
Hence up to taking a further subsequence we may assume
that $h_i(z)$ converges to some point in $\si$. Since
$\tau(h_i(z),h_i(u_i))
\rightarrow +\infty$ and $h_i(u_i) \rightarrow u$, lemma \ref{conv}
implies that $h_i(z) \rightarrow \eta_-(u)$.
In addition since $h_i(z) \cdot [0,+\infty)$
contains $h_i(u_i) \rightarrow u \in \mi$ 
and $h_i(u_i) = h_i(z) \cdot r_i$ with $r_i \rightarrow +\infty$,
it follows that 
every point in $u \cdot \rrrr$ is obtained as a limit
of points $q_i \in h_i(z) \cdot [0,+\infty)$.
Since $\pi^{-1}(P')$ is a closed set in $\mi$, 
this shows that $\pi(u \cdot \rrrr)$ 
is contained in $P'$ and therefore $u \cdot \rrrr$ is a 
$K$-quasigeodesic.
In the same way $h_i(y) \rightarrow \eta_+(v)$ and 
$v \cdot \rrrr$ is obtained as a limit of $h_i(y) \cdot (-\infty,0]$,
hence $v \cdot \rrrr$ is also a $K$-quasigeodesic,
see fig. \ref{dist}.

\begin{figure}
\centeredepsfbox{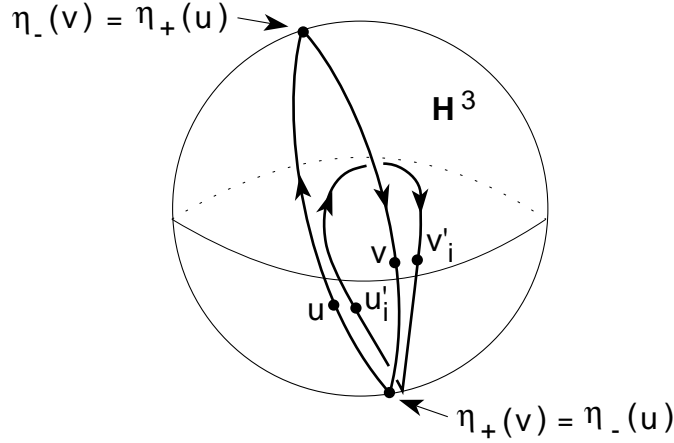}
\caption{Identification of ideal points in an orbit produces two
quasigeodesic orbits with same pair of ideal points.}
\label{dist}
\end{figure}

As $d(h_i(z), h_i(y))$ is a constant, it follows that
$\eta_+(v) = \eta_-(u)$.
Since both $v \cdot \rrrr$ and $u \cdot \rrrr$ are $K$-quasigeodesics,
this equality implies that $u$ and $v$ are not in the
same flow line of $\wwp$. 
But for each $i$, $h_i(v_i)$ and $h_i(u_i)$ are in the same
orbit of $\wwp$, $h_i(v_i) \rightarrow v$, and $h_i(u_i) \rightarrow
u$. This contradicts the Hausdorff orbit space condition.
We conclude that $\gamma_0$ and $\gamma_1$ do not have
the same ideal point.


\medskip
\noindent
{\underline {Case 2}} --- 
$x \cdot \reals$ intersects $\pi^{-1}(S)$ infinitely many times. 

Then there is a sequence $E_1 < \cdots < E_{J_0} \in {\mathcal E}$
such that $x \cdot \reals$ intersects each, and by
lemma \ref{LimitSetsDisjoint} 
we have $\ClBc(E_1) \intersect \ClFr(E_{J_0})
= \emptyset$. But $\eta_-(x) \in \ClBc(E_1)$ and $\eta_+(x) \in
\ClFr(E_{J_0})$ so $\eta_-(x) \ne \eta_+(x)$.
\end{proof}

Finally we prove property (c) of \ref{upgrade}.

\begin{proposition}{(continuity of extension)}{}
The map $\eta_+: \cover \Omega \rightarrow \si$
is continuous and similarly for $\eta_-$.
\label{cont}
\end{proposition}

\begin{proof}{}
We only consider $\eta_+$.

\smallskip
\noindent
{\underline {Case 1}} --- $x \cdot \RPlus$  intersects
$\pi^{-1}(S)$ infinitely often.

Let $E_1 < E_2 < \cdots$ be the elements of $\mathcal E$ that $x \cdot
\RPlus$ intersects. From lemma \ref{BoundedBelow},  $\Intersection_{i
\in {\bf N}} \ClFr(E_i)$ is a single point $p$ in $\bdy_\infinity
\cover P$. Since $\ClFr(E_i) \cap \bdy_\infinity\cover P$ is a
decreasing sequence of compact sets in $\bdy_\infinity\cover P$, then
for each $\delta > 0$, there is $i_0$ so that $\forall i > i_0$, the
diameter of $\ClFr(E_i) \cap \bdy_\infinity\cover P < \delta$. If $z$
is sufficiently near $x$ then $z \cdot \RPlus$ will intersect
$E_{i_0}$ therefore $\eta_+(z) \in \ClFr(E_i) \cap \bdy_\infinity\cover
P$. Hence the distance from $\eta_+(z)$ to $\eta_+(x)$ $< \delta$.
This implies that $\eta_+$ is continuous at $x$.

\smallskip
\noindent
{\underline {Case 2}} --- $x \cdot \RPlus$ has finite intersection
with $\pi^{-1}(S)$.

There is $s_1 \geq 0$ with $x \cdot [s_1,+\infty) \cap
\pi^\inverse(S) =
\emptyset$. Let $x' = x \cdot s_1$. Since 

$$x_i \rightarrow x  \ \ \Longleftrightarrow  \ \ 
x_i \cdot s_1 \rightarrow x',$$  

\noindent we may assume that $x \cdot \Rplus \cap \pi^{-1}(S)  =
\emptyset$. Consider a sequence $x_i \to x$, and let
$$ m_i = \bigl| x_i \cdot \Rplus  \intersect \pi^\inverse(S) \bigr|
$$
where $\bigl| V \bigr|$ denotes cardinality of $V$. Passing to a
subsequence, either $m_i$ takes a constant value $\le J_0$ for all
$i$, or $m_i > J_0$ for all $i$.

\smallskip
\noindent
{\underline {Case 2.1}} --- $m_i = 0$ for all $i$.

We may assume all $\gamma_i = x_i \cdot \RPlus$ and $\gamma_0 =
x \cdot \RPlus$ are in the $\delta_1$-neighborhood of $\Omega$. By
lemma \ref{uniform}, $\gamma_0$ and each $\gamma_i$ are
$K$-quasigeodesic. The proof in this case is finished by lemma
\ref{conv}.

Before addressing the case $m_i > 0$ we  need the following lemma.

\begin{lemma}{}{}
Suppose $m_i \geq j-1$ for all $i$. Let $v_i$ be the $j$-th
intersection of $x_i~\cdot \RPlus$ with $\pi^{-1}(S)$, if it exists,
otherwise let $v_i = \eta_+(x_i)$. Then $v_i$ converges to 
$\eta_+(x)$ in $\widetilde P \cup \partial _{\infty} \widetilde P$.
\label{inter}
\end{lemma}

\begin{proof}{} First suppose $j=1$, and let $v_i = x_i \cdot t_i$
with $t_i \in [0,+\infinity]$. By lemma \ref{uniform} the segments
or rays $[x_i,v_i)$ are all $K$-quasigeodesics. In addition
$t_i\rightarrow +\infty$, since $x \cdot \RPlus \cap \pi^{-1}(S) =
\emptyset$ and $x_i \rightarrow x$. By lemma \ref{conv} we conclude
that $z_i \rightarrow \eta_+(x)$.

\begin{figure}
\centeredepsfbox{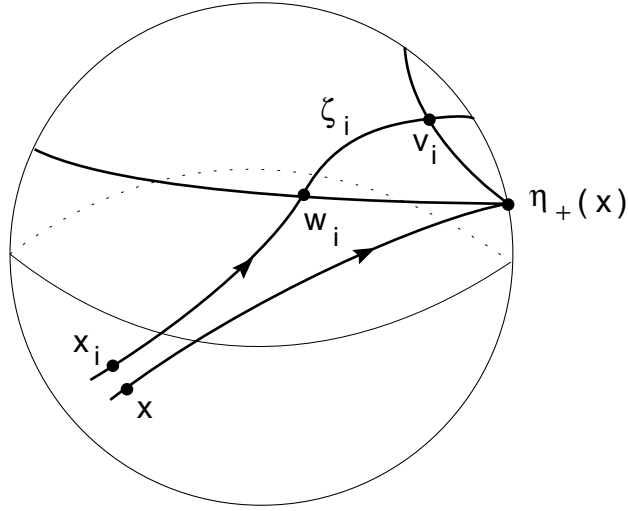}
\caption{Successive intersections with ${\mathcal E}$.}
\label{induction}
\end{figure}

Suppose now that $m_i \ge j-1$ for all $i$, and let $w_i$ be the
$(j-1)$-th intersection. We assume by induction that $w_i \to
\eta_+(x)$. Let $\zeta_i = [w_i,v_i]$, so each $\zeta_i$ is a
$K$-quasigeodesic, see fig. \ref{induction}.
If $\zeta_i$ is not escaping compact sets in $\cover
P$, then passing to a subsequence $\zeta_i$ accumulates on some $u \in
\cover P$, leading to a contradiction as follows: if $u\in
x\cdot\rrrr$, then since each $\zeta_i - \{ w_i \}$ is separated
from $x$ by at least one  element of ${\mathcal E}$, this would
imply that $x \cdot \rrrr_+$ has to intersect ${\mathcal E}$ contradiction.
In the other case,
$x$ and $u$ are not on the same flow line, contradicting the
Hausdorff orbit space condition.

The paths $\zeta_i$ therefore escape to infinity in $\cover P$. These
paths are $K$-quasigeodesics,
so their diameters are shrinking to zero
in a topological metric on 
$\widetilde P \cup \partial _{\infty} \widetilde P$.
Since $w_i \rightarrow
\eta_+(x)$ it then follows that $v_i \rightarrow \eta_+(x)$.
Induction now completes the proof of the lemma.
\end{proof}

\smallskip
\noindent
{\underline {Case 2.2}} --- $m_i = m$ is a constant with $0 < m \le
J_0$.

Applying lemma \ref{inter} with $j=m-1$ it follows that $\eta_+(x_i)
\to \eta_+(x)$.

\vskip .15in
\noindent
{\underline {Case 2.3}} --- $m_i > J_0$ for every $i$.

Let $y_i$ be the corresponding $(J_0 + 1)$-th intersection. Then $y_i
\rightarrow \eta_+(x)$. Let $C_i \in {\mathcal E}$ be the lift of $S$
containing $y_i$.

If there is an infinite subsequence $i_l$ with $C_{i_l} = C$ for all
$l$, then because $z_{i_l} \rightarrow \eta_+(x)$, it follows that
$\eta_+(x) \in \bdy_{\infty} C$. But for any $i$, the flow line through
$x_i$ intersects a sequence $E_1 < \cdots < E_{J_0} \in \mathcal E$ before
intersecting $C_i$. Since $x \cdot \Rplus \subset \sigma^b E_1$ then
$\eta_+(x) \in \sigma^b E_1$, but also $\eta_+(x) \in \bdy_\infinity
C \subset \sigma^f E_{J_0}$, and therefore $\sigma^b E_1
\intersect \sigma^f E_{J_0} \ne \emptyset$ contradicting lemma
\ref{LimitSetsDisjoint}. We may therefore assume up to taking a
further subsequence that all $C_i$ are distinct from each other.

Note that there is $E_0 \in \mathcal E$ such that $E_0 < C_i$ for all
$i$; take any $E_0$ so that $x \in \Fr(E_0)$. Applying lemma
\ref{BoundedBelow}, there is a subsequence $i_l$ so that 
that $\sigma^f C_{i_l}$ approaches a single point of
$\bdy_\infinity P$ in the Hausdorff topology. This point must be
$\eta_+(x)$, since $z_{i_l} \in \sigma^f(C_{i_l})$ and $z_{i_l} \to \eta_+(x)$.
Therefore the sequence $\eta_+(x_{i_l}) \in \sigma^f(C_{i_l})$ approaches
$\eta_+(x)$.

\smallskip
The  arguments of cases 2.1, 2.2 and 2.3,
show that given any sequence $y_i \rightarrow x$,
there is always a  subsequence $y_{i_l}$ for which $\eta_+(y_{i_l})
\rightarrow \eta_+(x)$. This implies that the original sequence
$\eta_+(y_i) \rightarrow \eta_+(x)$. This finally finishes the proof of
proposition \ref{cont}.
\end{proof}

Propositions \ref{ext}, \ref{distinct} and \ref{cont}
show that orbits in $\Omega_{n+1}$
satisfy conditions (a,b,c) of theorem \ref{upgrade}
as orbits in $\cover P \cup \partial _{\infty} \widetilde P$.
Since $\cover P \cup \partial _{\infty} \widetilde P$ embeds
continuously in $\hhh \cup \si$, it follows that they also satisfy
these conditions as seen in $\hhh \cup \si$.
Theorem \ref{upgrade} then implies that orbits of $\wwp$ in $\widetilde P$
are uniformly quasigeodesic in $\hhh$ (hence also in $\widetilde P$).
Induction now finishes the proof
of theorem $A$.

\medskip
\noindent
{\bf Remark:} \ 
In this article we use the inductive step of this section to
study
flows transverse to Reebless finite depth foliations.
By Novikov's theorem \cite{No}, $\Phi$ satisfies the 
$\epsilon, T$-cycles condition,
hence lemma \ref{weakqg}  can be applied, 
simplifying the proof of case 1 in proposition \ref{distinct}.
However, the inductive 
proof we give in this section has more general hypothesis, namely:
(1) $\wwp$ has Hausdorff orbit space, (2) $\Phi$ is well adapted to
a partial sutured manifold hierarchy, that is,
the smallest sutured manifold in the hierarchy, call it $P$ 
may not be a product sutured manifold; and (3) orbits entirely
contained in $P$ are uniform quasigeodesics.
There may be flows satisfying these properties,
which are not transverse to finite depth foliations.


\section{Quasigeodesic pseudo-Anosov flows}

Pseudo-Anosov flows are a generalization of suspension flows of
pseudo-Anosov surface homeomorphisms. These flows behave much like
Anosov flows, but they have finitely many singular orbits with
a prescribed behavior. In order to define pseudo-Anosov flows, first
we recall singularities of pseudo-Anosov surface homeomorphisms. 

Given $n \ge 2$, the quadratic differential $z^{n-2} dz^2$ on the complex
plane $\complex$ (see \cite{St} for quadratic differentials)
has a horizontal singular foliation $f^u$ with
transverse measure $\mu^u$, and a vertical singular foliation $f^s$
with transverse measure $\mu^s$. These foliations have $n$-pronged
singularities at the origin, and are regular and transverse to each
other at every other point of $\complex$. Given $\lambda > 1$, there is
a map $\psi \from \complex \to \complex$ which takes $f^u$ and $f^s$ to
themselves, preserving the singular leaves, stretching the leaves of
$f^u$ and compressing the leaves of $f^s$ by the factor $\lambda$. Let
$R_\theta$ be the homeomorphism $z \to e^{2 \pi \theta} z$ of
$\complex$. If $0 \le k < n$ the map $R_{k/n} \composed
\psi$ has a unique fixed point at the origin, and this defines the
local model for a {\em pseudohyperbolic fixed point}, with $n$-prongs
and with rotation $k$. Let $d_\Euclidean$ be the singular
Euclidean metric on
$\complex$ associated to the quadratic differential $z^{n-2} dz^2$, given by 
$$d_\Euclidean^2 = \mu_u^2 + \mu_s^2
$$ 
Note that 
$$(R_{k/n} \composed \psi)^* d_\Euclidean^2 = \lambda^{-2} \mu_u^2 +
\lambda^2 \mu_s^2
$$

Now consider the mapping torus $N = \complex \cross \reals / (z,r+1)
\sim (R_{k/n} \composed \psi(z),r)$, with suspension flow $\Psi$
arising from the flow in the $\reals$ direction on $\complex \cross
\reals$. The suspension of the origin defines a periodic orbit
$\gamma$ in $N$, and we say that $(N,\gamma)$ is the local model for a
{\em pseudohyperbolic periodic orbit}, with $n$ prongs and with
rotation $k$. The suspension of the foliations $f^u,f^s$ define
2-dimensional foliations on $N$, singular along $\gamma$, called the
{\em local weak unstable and stable foliations}. Note that there is a
singular Riemannian metric
$ds$ on $\complex \cross \reals$ that is preserved by the gluing
homeomorphism $(z,r+1) \sim (R_{k/n} \composed \psi(z),r)$, given by
the formula
$$ ds^2 = \lambda^{-2t} \mu_u^2 + \lambda^{2t} \mu_s^2 + dt^s
$$ 

The metric $ds$ descends to a metric on $N$ denoted $ds_N$.

Let $\Phi$ be a flow on a closed, oriented 3-manifold $M$. We say
that $\Phi$ is a {\em pseudo-Anosov flow} if the following are
satisfied:

\begin{itemize}

\item For each $x \in M$, the flow line $t \to \Phi(x,t)$ is $C^1$,
and the tangent vector bundle $D_t \Phi$ is $C^0$.

\item There is a finite number of periodic orbits $\{ \gamma_i \}$,
called {\em singular orbits}, such that the flow is smooth off of the
singular orbits.

\item Each singular orbit $\gamma_i$ is locally modelled on a
pseudo-hyperbolic periodic orbit. More precisely, there exist $n,k$
with $n \ge 3$ and $0 \le k < n$, such that if $(N,\gamma)$ is the
local model for an pseudo-hyperbolic periodic orbit with $n$ prongs
and with rotation $k$, then there are neighborhoods $U$ of $\gamma$
in $N$ and $U_i$ of $\gamma_i$ in $M$, and a diffeomorphism $f
\from U \to U_i$, such that $f$ takes orbits of the semiflow $R_{k/n} \composed
\psi \restrictedto U$ to orbits of $\Phi \restrictedto U_i$.

\item There exists a path metric $d_M$ on $M$, such that $d_M$ is a
smooth Riemannian metric off of the singular orbits, and for a
neighborhood $U_i$ of a singular orbit $\gamma_i$ as above, the
derivative of the map $f \from (U - \gamma) \to (U_i - \gamma_i)$ has
bounded norm, where the norm is measured using the metrics $ds _N$ on
$U$ and $d_M$ on $U_i$.

\item On $M - \Union \gamma_i$, there is a continuous splitting of the
tangent bundle into three 1-dimensional line bundles $E^u \oplus E^s
\oplus T \Phi$, each invariant under $\Phi$, such that $T \Phi$ is
tangent to flow lines, and for some constants 
$\nu > 1, \theta >  1$ we have

\begin{enumerate}

\item If $v \in E^u$ then $|D\Phi_t(v)| \le 
\theta \nu^t |v|$ for $t<0$
\item If $v \in E^s$ then $|D\Phi_t(v)| \leq \theta \nu^{-t} |v|$ for
$t>0$

\end{enumerate}

{\noindent
where norms of tangent vectors are measured using the metric $d_M$.}

\end{itemize}

With the definition formulated in this manner, the entire theory of
Anosov flows can be mimicked for pseudo-Anosov flows. In particular,
a pseudo-Anosov flow $\Phi$ has a 2-dimensional weak unstable
foliation $\fu$ tangent to $E^u \oplus T \Phi$ away from the singular
orbits, and a 2-dimensional weak stable foliation $\fs$ tangent to
$E^s \oplus T \Phi$. These foliations are singular along the singular
orbits of $\Phi$, and regular everywhere else. In the
neighborhood $U_i$ of an $n$-pronged singular orbit $\gamma_i$, the
images of $\fs$ and $\fu$ in the model manifold $N$ are identical with
the local weak stable and unstable foliations. 

The restriction of
$\fs$ to $M - $ (singular orbits) defines a true foliation; a
complete leaf of this foliation is called a {\em nonsingular leaf} of
$\fs$; an incomplete leaf may be completed by adding a singular orbit
$\gamma$ of $\Phi$, and the result is called a {\em singular leaf} of
$\fs$ abutting $\gamma$. Singular and nonsingular leaves of $\fu$ are
similarly defined. The bare term ``leaf'' means either a nonsingular
or a singular leaf. For some small neighborhood $W$ of any point $x$
lying on a singular orbit $\gamma$, the singular leaves of $\fs$
divide $W$ into $n$ {\em sectors}, each sector parameterized by $D^2
\cross I$, so that the restriction of $\fs$ to the sector agrees with
the foliation by level discs $D^2 \cross t$; and similarly for $\fu$.
All the terms defined here apply as well to the lifted singular
foliations $\fns, \fnu$ in $\mi$.

\begin{proposition}{}{}
If $\Phi$ is a pseudo-Anosov flow in $M^3$ then the
orbit space $\oo$ of $\wwp$ is homeomorphic to $\rrrr^2$.
\end{proposition}

\begin{proof}{}
Let $\fs, \fu$ be the singular stable and unstable foliations of
$\Phi$. A short embedded path $\alpha$ is a {\em quasi-transversal} of
$\fs$ if one of the following happens: either $\alpha$ is transverse
to $\fs$; or there is a singular orbit $\gamma$ and $x \in
\interior(\alpha) \intersect \gamma$ such that the closure of each
half of $\alpha - x$ either is transverse to $\fs$ or lies on a
singular leaf of $\fs$, and the two halves do not both lie in the
closure of any sector of $\fs$ at $x$. Quasitransversals of $\fu$ are
similarly defined. We need the following facts about $\Phi$, $\fs$ and
$\fu$.

\begin{enumerate}

\item Each periodic orbit of $\Phi$ is homotopically nontrivial.

\item A quasitransversal $\rho$ of $\fu$ is not path homotopic into a
leaf of $\fu$; and similary for $\fs$.

\item The inclusion map of every leaf of $\fns$ or $\fnu$ into $\mi$
is proper.

\end{enumerate}

These facts are proved using the theory of essential laminations
\cite{Ga-Oe}; we give the argument for $\fs$. Split $\fs$ along
singular leaves to produce a lamination $\ls$. There is one
complementary component $T_\gamma$ of $\ls$ for each singular orbit
$\gamma$ of $\Phi$, and if $\widehat T_\gamma$ is the metric completion of
$T_\gamma$ then $\widehat T_\gamma$ is homeomorphic to a solid torus with
an $(n,k)$ torus knot removed from the boundary, where $\gamma$ has $n$
prongs and rotation $k$. Since $n \ge 3$ for all $\gamma$ it follows
that $\ls$ is an essential lamination. There is a map $f \from M \to
M$ homotopic to the identity taking $\ls$ onto $\fol^s$,
taking each leaf of $\ls$ not
on the boundary of some $T_\gamma$ homeomorphically onto a
non-singular leaf of $\fs$, and taking each leaf on the boundary of
each $T_\gamma$ onto a union of singular leaves abutting
$\gamma$. Given a leaf $L \subset \bdy T_\gamma$ if gcd$(n,k) \ne
1$ then $L$ is mapped onto a union of two singular leaves and the map
is 1-1 off of $\gamma$; and if gcd$(n,k) = 1$ then $L$ maps onto one
singular leaf and the map is 2-1 off of $\gamma$; in either case,
$f^\inverse(\gamma) \intersect L$ is a core curve of the annulus $L$.

Properties (1--3) follow from results of \cite{Ga-Oe}. Property (1)
follows because each periodic orbit of $\Phi$ is homotopically
nontrivial in a leaf of $\fs$, and because leaves of $\ls$ are
$\pi_1$-injective. Property (2) follows because if $\rho'$ is a path
transverse to $\ls$, and if the closure of each component of $\rho' -
\ls$ is not path homotopic into a leaf of
$\ls$, then $\rho'$ is not path homotopic into a leaf of $\ls$.
Property (3) follows because leaves of $\lns$ and $\lnu$ include
properly in $\mi$.

In order to prove the proposition,
the key facts to prove are that $\oo$ is a two dimensional manifold and
that it is Hausdorff.

Suppose there are $\epsilon, T$-cycles of $\Phi$ for $\epsilon$
arbitrarily small and $T$ arbitrarily big. We can then choose $x_i
\rightarrow x \in \mi$ and $t_i \to +\infinity$ so that $y_i = x_i
\cdot t_i \to x$. Let $\alpha_i = x_i \cdot [0,t_i]$. 
If $x_i, y_i$ are in the same local orbit of $\wwp$ near $x$,
then $\alpha_i$ can be easily completed to a closed orbit,
contradicting  property (1) above. 
If the endpoints of $\alpha_i$ are not in the same local
orbit and $\epsilon$ is sufficiently small, 
then the endpoints of $\alpha_i$ may be joined by a path $\beta$
which is a quasi-transversal to either $\fns$ or $\fnu$; but
$\beta$ is path homotopic to $\alpha_i$, contradicting property (2)
above.



This implies that $\oo$ is locally $2$-dimensional. We next prove that
it is Hausdorff. The proof given for Anosov flows in \cite{Fe1} goes
through almost verbatim, with slight changes to take pseudo-Anosov
behavior into account, or we may proceed as follows.

Let $x_i \rightarrow x$, $y_i = x_i \cdot t_i \rightarrow y$ in
$\mi$.  We first want to show that $t_i$ is bounded. Assume then
up to taking a subsequence that $t_i \rightarrow +\infty$. If $x$ lies
on a singular orbit $\gamma$ of $\Phi$, pass to a subsequence so that
$x_i$ lies in a single sector of $\fns$ near $x$. Now vary $x_i$
to nearby $z_i$, so that $z_i \rightarrow x$, $z_i \in \ws(x_i) \cap
\wu(x)$. Since flow lines in the stable foliation converge together
exponentially it follows that there are $s_i \rightarrow +\infty$ with
$w_i = z_i \cdot s_i \rightarrow y$. The sequence $z_i$ forms a bounded
subset of a leaf $L$ of $\wu(x)$, and hence the sequence $z_i \cdot
s_i$ leaves every compact subset of $L$,
because flow lines are properly embedded in
the leaf $L$ --- the flow induces a product
structure  in $L$. But by (3) above the inclusion
map $L \inject \mi$ is proper, contradicting that $z_i \cdot s_i$
converges in $\mi$. This contradiction shows that the $t_i$ are
bounded.

As the $t_i$ are bounded, we may assume up to subsequence
that $t_i \rightarrow t_0$. Then $y_i \rightarrow x \cdot t_0 = y$.
This shows that $\oo$ is Hausdorff and it follows that
it is homeomorphic to $\rrrr^2$.
\end{proof}

Given any closed, irreducible,  orientable $3$-manifold $M$
with  $H_2(M) \not = 0$, Gabai \cite{Ga1} constructed many
finite depth foliations in $M$: for any $\zeta \not = 0$ in
$H_2(M)$ there is a taut, finite depth foliation $\fol$
in $M$ so that the compact leaves of $\fol$ represent $\zeta$.

If in addition $M$ is atoroidal, then Mosher \cite{Mo4} constructed
pseudo-Anosov flows in $M$ which are {\em almost tranverse} to $\fol$.
These flows become transverse to $\fol$ after blow up of finitely
many singular orbits $\gamma_i, 1 \leq i \leq i_0$ of $\Phi$.
The blown up flow is denoted by $\Phi^b$. The blow up transforms the
orbit $\gamma_i$ into  $Z_i = (T_i \times [0,1])/f$, where $T_i$ is a
finite simplicial tree, and $f$ is a homeomorphism of $T_i \times \{ 1
\}$ to $T_i \times \{ 0 \}$ which sends edges to edges.
Each edge $E$ of $T_i$ eventually returns to itself, producing an
annulus $A$ in $Z_i$. The flow in $A$ is as follows: the boundary
consists of two orbits coherently oriented, and the interior orbits
spiral from one boundary circle to the other without forming a two
dimensional Reeb component. Clearly there is a global section to
$\Phi^b$ restricted to $Z_i$. The flow $\Phi^b$ is transverse to
$\fol$.

In addition $\Phi^b$ is semiconjugate to $\Phi$: there is a continuous
map $q \from M \to M$, so that 

\begin{itemize}
\item $q$ is homotopic to the identity,
\item 
$q$ takes orbits of $\Phi^b$ to orbits of $\Phi$ preserving
orientation,
\item
$q$ is $C^1$ along orbits of $\Phi^b$,
\item
$q$ is one to one
except in $\cup_{1 \leq i \leq i_0} Z_i$,
\item $q(Z_i) = \gamma_i$.
\end{itemize}

\begin{proposition}{}{}
The orbit space $\oo^b$ of $\wwp^b$ is homeomorphic to $\rrrr^2$.
\label{orbit}
\end{proposition}

\begin{proof}{}
Since $\Phi^b$ is transverse to $\fol$, then if $\epsilon>0$ is very
small and $T>0$ very large, it follows that an $\epsilon, T$ cycle
can be perturbed to be transverse to $\fol$, hence it is not null
homotopic by Novikov's theorem \cite{No}. 
This shows that $\oo^b$ is locally two
dimensional. We now show that $\oo^b$ is Hausdorff.

Lift $q$ to a map $\widetilde q \from \mi \rightarrow \mi$,
which is boundedly homotopic to the identity map. Let $\tau(x,y)$ be
the difference in flow parameter between $x,y$ on the same flow line
of $\wwp$, and let $\tau^b(x,y)$ be similarly defined for $\wwp^b$.

\smallskip
\noindent
{\underline {Claim}} ---
There are $m_0, m_1 > 0$ so that if $x, y$ are in the same flow line
of $\wwp^b$, then  $m_0 \tau^b(x,y) < \tau(\tilde q(x), \tilde q(y)) <
m_1 \tau^b(x,y)$.

Given $z \in M$, let $g(z)$ be the derivative at $z$ of $q$ restricted
to the flow line through $z$. Then since $g(z)$ is continuous and
positive and $M$ compact, there are $m_0$ and $m_1 > 0$ which are the
maxima and minima of $g$ in $M$. The claim follows.

%
%
%
%
\smallskip

Let now $x_i \in \mi, y_i = \wwp^b_{t_i}(x_i)$ so that $y_i
\rightarrow y$ and $x_i \rightarrow x$. Then $\tilde q(x_i) \rightarrow
\tilde q(x)$ and $\tilde q(y_i) \rightarrow \tilde q(y)$.  Since $\Phi$
is pseudo-Anosov, the previous proposition shows that 
the orbit space $\oo$ of $\wwp$ is
Hausdorff, hence $\tilde q(y) = \wwp_s(\tilde q(x))$. Since $\oo$ is
homeomorphic to $\rrrr^2$, it follows that $\tilde q(y_i) = 
\wwp_{s_i}(\tilde q(x_i))$ where $s_i \rightarrow s$.

Partitioning $s_i$ into a uniformly bounded number of subsegments of
length $\leq 1$, and using the claim above, it follows that $y_i =
\wwp^b_{r_i}(x_i)$ where $r_i$ is uniformly bounded. 
Up to taking a 
subsequence, assume that $r_i$ converges to some $r$, and hence
$y_i\rightarrow\wwp^b_r(x)$. Consequently
$y = \wwp^b_r(x)$. This shows that
$\oo^b$ is Hausdorff. It follows that $\oo^b$ is homeomorphic to
$\rrrr^2$.
\end{proof}

Now we prove the main theorem:

\medskip
\noindent{\bf Main theorem} \ 
{\em Let $M^3$ be closed, orientable, irreducible with non zero
second betti number. Let $\zeta \not = 0 \in H_2(M)$ and 
$\fol$ be a taut finite depth foliation with compact
leaves representing $\zeta$. Let $\Phi$ be a pseudo-Anosov flow
which is almost transverse to $\fol$. Then $\Phi$ is quasigeodesic.}
\smallskip

\begin{proof}{} 
By the previous proposition the flow $\widetilde\Phi^b$ on $\widetilde M$
has Hausdorff orbit space. 

Suppose first that $\fol$ has a compact leaf which is not a fiber.
Then $\Phi^b$ is quasigeodesic by theorem A. The semiconjugacy $\tilde
q$ from $\widetilde \Phi^b$ to $\widetilde \Phi$ moves points a
uniformly bounded distance in $\mi$. Since there is also a bound on
how much the map $\tilde q$ expands the lengths of flow lines (given
by the claim in the previous proposition), it follows that flow lines
of $\wwp$ are uniform quasigeodesics.

Now suppose that each compact leaf of $\fol$ is a fiber. If $\fol$ is
a fibration over the circle, then $\Phi^b$ is quasigeodesic by
\cite{Ze}, and the arguments of the previous paragraph show that
$\Phi$ is quasigeodesic. 

Suppose $\fol$ is not a fibration over the circle. We claim that
$\fol$ can be replaced by another finite depth foliation of smaller
depth that is still transverse to $\Phi^b$; continuing by induction,
eventually $\Phi^b$ is transverse to a fibration over the circle, and
we are done.

We sketch a proof of the claim. Let $\fol$ have depth $n$. Note that
each component of $\fol - \fol^0$ is homeomorphic to the product of a
depth 0 leaf crossed with an interval; similarly each component of
$\fol - \fol^k$ is homeomorphic to the product of a depth $k$ leaf
crossed with an interval. Note also that the restriction of $\fol$ to a
component of $\fol - \fol^{n-1}$ is a fibration over the circle with
fiber a depth $n$ leaf $L$. It follows that if $f
\from L \to L$ is the first return map of $\Phi^b$, then there is a
{\em translation map} $g \from L \to L$, i.e.\ a map which generates a
properly discontinuous, free action of
$\integers$, such that $f$ is isotopic to $g$ by a compactly supported
isotopy. Any compact subset of $L$
invariant under $f$ is contained in the support of the isotopy from
$f$ to $g$, and hence there is a maximal compact invariant set $C$ of
$f$. Using the fact that $\Phi^b$ is a blown up pseudo-Anosov flow,
the set $C$ is nonempty if and only if $f$ has periodic points.
Moreover, any periodic points are non-removable in the proper
isotopy class of $f$; the proof given for pseudo-Anosov surface
homeomorphisms in \cite{Bi-Ki} works just as well in the present
context. Since $g$ has no periodic points, it follows that $C =
\emptyset$, and therefore $f$ is itself a translation map. The
component of $\fol - \fol^{n-1}$ containing $L$ may therefore be
refoliated by leaves of depth $n-1$, staying transverse to $\Phi^b$.
Doing this for each component of $\fol - \fol^{n-1}$ proves the claim.
\end{proof}

{\footnotesize
{

\setlength{\baselineskip}{0.05cm}

}
}

\end{document}